\documentclass[11pt, reqno]{amsart}
\usepackage{amsmath}
\usepackage{amssymb}

\numberwithin{equation}{section}

\newtheorem{theorem}{Theorem}[section]
\newtheorem{lemma}[theorem]{Lemma}
\newtheorem{proposition}[theorem]{Proposition}
\newtheorem{corollary}[theorem]{Corollary}
\newtheorem{definition}[theorem]{Definition}
\newtheorem*{dbp}{The Askey-Wilson bispectral problem}

\theoremstyle{definition}
\newtheorem{remark}[theorem]{Remark}

\newcommand{\thref}[1]{Theorem \ref{#1}}
\newcommand{\leref}[1]{Lemma \ref{#1}}
\newcommand{\prref}[1]{Proposition \ref{#1}}
\newcommand{\reref}[1]{Remark \ref{#1}}
\newcommand{\deref}[1]{Definition \ref{#1}}

\newcommand{\coref}[1]{Corollary \ref{#1}}
\newcommand{\seref}[1]{Section \ref{#1}}

\newcommand{\phift}[8]{{}_4\phi_3\left[\begin{matrix} #1 , #2, #3, #4 \\
#5, #6, #7 \end{matrix}\,; q,#8\right]}
\newcommand{\Wes}[7]{{}_8W_7\left(#1; #2, #3, #4, #5, #6; q, #7\right)}

\newcommand{\RE}{{\mathcal R}\{q^{\gamma},E\}}

\newcommand{\Lg}{L_{a,b,c,d}(\gamma,E)}
\newcommand{\Bz}{B_{a,b,c,d}(z,D_z)}

\newcommand{\La}[1]{L_{#1,b,c,d}}
\newcommand{\Pa}{P_{a,b,c,d}}
\newcommand{\Qa}{Q_{a,b,c,d}}
\newcommand{\Ag}[1]{A_{#1;\,a,b,c,d}}
\newcommand{\Bg}{B_{a,b,c,d}}
\newcommand{\Cg}[1]{C_{#1;\,a,b,c,d}}
\newcommand{\Aa}[1]{A_{\gamma;\,#1,b,c,d}}

\newcommand{\Ca}[1]{C_{\gamma;\,#1,b,c,d}}
\newcommand{\Aza}[1]{A_{a,b,c,d}(z)}
\newcommand{\Az}{A_{a,b,c,d}(z)}
\newcommand{\rga}{r_{\gamma}(a,b,c,d;z)}
\newcommand{\ra}[1]{r_{\gamma}(#1,b,c,d;z)}
\newcommand{\raz}[2]{r_{\gamma}(#1,b,c,d;#2)}

\newcommand{\rgaf}[2]{{\mathfrak{r}}_{#1}(#2,b,c,d;z)}
\newcommand{\sga}{s_{\gamma}(a,b,c,d;z)}
\newcommand{\sa}[1]{s_{\gamma}(#1,b,c,d;z)}
\newcommand{\saz}[2]{s_{\gamma}(#1,b,c,d;#2)}

\newcommand{\sgaf}[2]{{\mathfrak{s}}_{#1}(#2,b,c,d;z)}
\newcommand{\at}{{\tilde{a}}}
\newcommand{\bt}{{\tilde{b}}}
\newcommand{\ct}{{\tilde{c}}}
\newcommand{\dt}{{\tilde{d}}}
\newcommand{\zt}{{\tilde{z}}}
\newcommand{\gt}{{\tilde{\gamma}}}
\newcommand{\Fg}[1]{F_{#1}(a,b,c,d;z)}
\newcommand{\Fz}[1]{F_{\gamma}(a,b,c,d;#1)}
\newcommand{\Fgt}[1]{F_{#1}(\at,\bt,\ct,\dt;\zt)}

\newcommand{\cB}{{\mathcal{B}}}
\newcommand{\cK}{{\mathcal{K}}}
\newcommand{\cL}{{\mathcal{L}}}
\newcommand{\cP}{{\mathcal{P}}}
\newcommand{\cQ}{{\mathcal{Q}}}
\newcommand{\cT}{{\mathcal{T}}}

\newcommand{\tQ}{{\tilde{\mathcal Q}}}

\newcommand{\hL}{{\hat{L}}}
\newcommand{\hB}{{\hat{B}}(z,D_z)}
\newcommand{\hPsi}{{\hat{\Psi}}(\gamma,z)}
\newcommand{\hLambda}{{\hat{\Lambda}}_{\gamma}}

\newcommand{\bm}{{\bar{m}}}
\newcommand{\fa}{{\mathfrak{a}}}
\newcommand{\fb}{{\mathfrak{b}}}
\newcommand{\fc}{{\mathfrak{c}}}
\newcommand{\fd}{{\mathfrak{d}}}
\newcommand{\ff}{{\mathfrak{f}}}
\newcommand{\fg}{{\mathfrak{g}}}
\newcommand{\fh}{{\mathfrak{h}}}
\newcommand{\fL}{{\mathfrak{L}}}
\newcommand{\fA}{{\mathfrak{A}}_{\gamma}}
\newcommand{\fB}{{\mathfrak{B}}_{\gamma}}
\newcommand{\fC}{{\mathfrak{C}}_{\gamma}}

\newcommand{\lL}{\lambda_{\Lambda_{\gamma}}}
\newcommand{\rL}{\rho_{\Lambda_{\gamma}}}
\newcommand{\dI}{{\Delta}^{I}}

\newcommand{\Id}{{\mathrm{Id}}}
\newcommand{\Ker}{{\mathrm{Ker}}\,}
\newcommand{\Span}{{\mathrm{Span}}}

\newcommand{\Z}{\mathbb{Z}}
\newcommand{\N}{\mathbb{N}}

\begin{document}
\title{Askey-Wilson type functions, with bound states}

\begin{abstract}
The two linearly independent solutions of the three-term
recurrence relation of the associated Askey-Wilson polynomials,
found by Ismail and Rahman in \cite{IR}, are slightly modified so
as to make it transparent that these functions satisfy a beautiful
symmetry property. It essentially means that the geometric and the
spectral parameters are interchangeable in these functions. We
call the resulting functions the Askey-Wilson functions. Then, we
show that by adding bound states (with arbitrary weights) at
specific points outside of the continuous spectrum of some
instances of the Askey-Wilson difference operator, we can generate
functions that satisfy
a doubly infinite three-term recursion relation and are also
eigenfunctions of $q$-difference operators of arbitrary orders.
Our result provides a discrete analogue of the solutions of
the purely differential version of the bispectral problem that
were discovered in the pioneering work \cite{DG} of Duistermaat and
Gr\"unbaum.
\end{abstract}

\author[L.~Haine]{Luc Haine}
\address{L.H., Department of Mathematics, Universit\'e catholique de
Louvain, Chemin du Cyclotron 2, 1348 Louvain-la-Neuve, Belgium}
\email{haine@math.ucl.ac.be}

\author[P.~Iliev]{Plamen Iliev}
\address{P.I., Department of Mathematics, University of California,
Berkeley, CA 94720-3840, USA}
\email{iliev@math.berkeley.edu}

\date{April 24, 2003}

\subjclass[2000]{33D45, 37K10, 14H70, 39A70, 39A13}

\maketitle


\section{Introduction}\label{se1}
Orthogonal polynomials which are eigenfunctions of
a differential operator have a long history. When
the differential operator is of order two, Bochner (1929) \cite{Bo} proved
that this
property characterizes the so-called classical orthogonal
polynomials, linked with the names of Hermite, Laguerre and Jacobi.
The general problem was raised by H.L. Krall \cite{Kr1} in 1938. He proved
that the differential operator has to be of even order and in
\cite{Kr2} he obtained a complete classification for the case of
an operator of order four.

In \cite{GH1}, it was observed that Krall's result can be reformulated in a
striking manner. The (semi-infinite) tridiagonal matrices which encode the
three-term recursion relation satisfied by Krall's polynomials,
are obtained by adding a bound state (with an arbitrary weight) at the boundary
of the continuous spectrum of the tridiagonal matrices associated
with some instances of the Laguerre and the Jacobi polynomials.
This way of looking at Krall's result leads to a very efficient construction
of his polynomials, by means of the so-called (discrete) Darboux
transformation, one of the basic tools in the theory of
solitons. Roughly, the method consists in factorizing the tridiagonal matrix
(appropriately shifted) as a product of a
lower and an upper matrix, producing then a new matrix by
permuting the factors. One of the nice outcomes of this approach
is that the Darboux process can be iterated, leading to orthogonal
sequences of polynomials satisfying differential equations,
with a moment functional given by a weight distribution
involving not only Dirac's delta function but also any of its
derivatives. We refer the reader to \cite{GHH} and \cite{GY} for
the precise formulation of the results in the context of the
Laguerre and the Jacobi polynomials respectively. For a sample of
other works related to Krall's problem or the role of the Darboux
transformation in the theory of orthogonal polynomials, the reader
can consult \cite{EKLW}, \cite{HI1}, \cite{HI3}, \cite{KK}, \cite{Ko}, \cite{L},
\cite{SZ} and \cite{Z}.

This paper deals with a discrete version of Krall's problem, replacing the
differential operator by a $q$-difference operator. The
polynomials that one meets here are of more recent vintage. The
celebrated Askey-Wilson polynomials \cite{AW} form the most general family
of orthogonal polynomials which are also eigenfunctions of a
second order $q$-difference operator, see \cite{AA}, \cite{GH2}.
In fact, we address an extension of the problem, by allowing for a \emph{doubly
infinite} three-term recursion relation, instead of a
semi-infinite one:
\begin{dbp} To determine all doubly infinite tridiagonal matrices
$L$, for which some family of eigenfunctions $\Psi(n,z)$
satisfying
\begin{equation}                        \label{1.1}
A_{n}\Psi(n+1,z)+B_{n}\Psi(n,z)+C_{n}\Psi(n-1,z)=(z+z^{-1})\Psi(n,z),
\end{equation}
is also a family of eigenfunctions of a $q$-difference operator $B$ in
the spectral variable $z$, i.e.
\begin{equation}                        \label{1.2}
\sum_{\text{finitely many $i$'s$\;\in\Z$}}D_{i}(z)\Psi(n,q^{i}z)
=\Lambda(n)\Psi(n,z).
\end{equation}
\end{dbp}

The motivation for allowing for doubly infinite tridiagonal matrices stems
from the pioneering work of Duistermaat and Gr\"unbaum \cite{DG},
where a purely differential version of the problem (with the
tridiagonal operator $L$ replaced by the Schr\"odinger operator)
was posed and completely solved.
The problem above was raised in \cite{GH3} by one of us in
collaboration with F.A. Gr\"unbaum, and solved there
in the case of a $q$-difference operator of order two. The
result is that, by allowing for doubly infinite
tridiagonal matrices, the only solutions are provided by what we
proposed to call the \emph{associated Askey-Wilson matrices}.
The entries of these matrices are obtained
by making an \emph{arbitrary} shift in the coefficients of the
recursion relations satisfied by the Askey-Wilson polynomials,
extending them over \emph{all} integers. Precisely,
\begin{equation}                        \label{1.3}
A_{n}=\tilde{A}_{n+\varepsilon},\quad
B_{n}=\tilde{B}_{n+\varepsilon},\quad
C_{n}=\tilde{C}_{n+\varepsilon},\quad n\in\mathbb{Z},
\end{equation}
where $\tilde{A}_{n},\tilde{B}_{n},\tilde{C}_{n}$, denote the
coefficients of the three-term recursion relation satisfied by the
Askey-Wilson polynomials.
The coefficients \eqref{1.3} define the so-called
\emph{associated Askey-Wilson polynomials} when the conditions
$\Psi(-1,z)=0$ and $\Psi(0,z)=1$ are imposed in \eqref{1.1}.
The associated Askey-Wilson polynomials were extensively studied
by Ismail and Rahman in the remarkable paper \cite{IR}.
However, as soon as $\varepsilon \neq 0$, the functions
$\Psi(\gamma,z),\gamma\equiv n+\varepsilon$, that solve \eqref{1.1} and
\eqref{1.2}, with a
second-order $q$-difference operator, are \emph{not} given by the
associated polynomials, but rather by a \emph{two-dimensional}
space of contiguous basic hypergeometric \emph{functions}, that
we called the \emph{Askey-Wilson functions}.
In \cite{GH3}, these functions were constructed
recursively in terms of an arbitrary solution of a second-order
$q$-difference equation (called there the Gauss-Askey-Wilson
equation), but explicit formulas were not obtained.

Our first result, \thref{th2.1}, is to give explicit formulas
for the Askey-Wilson functions in terms of basic hypergeometric series.
For this, we use two linearly independent eigenfunctions of the Askey-Wilson
second order difference operator $L\equiv\La{a}(\gamma,E)$ in terms of
very-well-poised $_8\phi_7$ functions, that were found by Ismail and Rahman
\cite{IR}, \cite{R2}. The notation $\La{a}(\gamma,E)$ reminds that the
operator $L$ depends on four arbitrary parameters (that are usually denoted by
$a,b,c,d$), $\gamma$ stands for $n+\varepsilon$ and $E$ denotes the
customary shift operator. We show that by multiplying the
Ismail-Rahman functions by an appropriate $z$ dependent
(but $\gamma$ independent) factor, these functions satisfy a beautiful
symmetry property, which we shall refer to as duality.
It essentially means that the geometric and the
spectral parameters $\gamma$ and $z$ of the eigenfunctions are
interchangeable, if we introduce an appropriate involution on the
parameters $a,b,c,d$. A special case of this duality was considered
by Koelink and Stokman \cite{KSt}, in their study of the
Askey-Wilson function transform. It is important to note that
these authors focus on a \emph{very specific} eigenfunction of the Askey-Wilson
difference operator, which they call the Askey-Wilson function,
because it is a meromorphic continuation of the Askey-Wilson
polynomial in its degree. The same function was studied previously
by Suslov \cite{Su1}, \cite{Su2}, who established  Fourier-Bessel
type orthogonality relations for it.

Sections 3 to 7 form the core of the paper. The main result can be found in
Section 6, \thref{th6.2} and \coref{co7.6}. We prove that starting with
\emph{some} instances of the associated Askey-Wilson matrices, and adding
bound states (with an arbitrary weight) at a collection of \emph{specific}
points \emph{outside of} the continuous spectrum of these
operators, we can construct solutions of the Askey-Wilson
bispectral problem, as stated at the beginning of this
introduction. We naturally call these solutions \emph{Askey-Wilson type
functions}. When $\varepsilon=0$ in \eqref{1.3}, we obtain in this way
orthogonal polynomials which are eigenfunctions of $q$-difference operators
of arbitrary orders, thus providing $q$-analogues of Krall's orthogonal
polynomials. As shown by formula \eqref{6.18} in \thref{6.2}
(see also \coref{co7.6}), when $q\to 1$, all the added bound states accumulate
at the boundary $\pm 2$ of the continuous spectrum $[-2,2]$ of the
Askey-Wilson second order difference operator. It is interesting to note that
most works on Krall's original problem deal with measures obtained by adding
delta functions (with arbitrary masses) to the measure of orthogonality of
some instances of the classical
orthogonal polynomials, at the boundary of the interval of
orthogonality, discarding the possibility of adding also derivatives of the
delta function. To the best of our knowledge, situations involving
derivatives of the delta function were first contemplated in
\cite{GHH}. In this case we don't have  an orthogonality \emph{measure}.
Our result shows that, in the $q$-difference case, the bound states
\emph{split up} and we do get an orthogonality
measure for the resulting polynomials. We shall now describe our strategy
to establish \thref{th6.2}, the proof of which is prepared by the material
exposed in Sections 3 to 5 of the paper.

The basic technique to establish \thref{th6.2} is the so-called
method of bispectral Darboux transformations which was developed
by Bakalov, Horozov and Yakimov \cite{BHY1}, \cite{BHY2},
and by Kasman and Rothstein \cite{KR}, in
relation with a program aiming at describing all bispectral commutative rings
of differential operators. In \cite{GHH}, the method was adapted to
attack more systematically Krall's original problem. Though it
was successful to produce extensions of the so-called
Krall-Laguerre polynomials, it led quickly to serious
computational difficulties, when applied starting with the Jacobi
polynomials. In the context of the Jacobi polynomials, the difficulties
were overcome by Gr\"unbaum and Yakimov \cite{GY}, by introducing
a new idea that we shall explain below.
Section 3 explains the method of bispectral Darboux transformations.
Roughly, the method consists in finding the most general
"rational" factorization of some constant coefficient polynomial in $\La{a}$
\begin{equation}                        \label{1.4}
\mathcal{L}=\fh(\La{a})=\prod_{i=1}^{m}(\La{a}-x_{i}\Id)=\mathcal{P}\mathcal{Q},
\end{equation}
for appropriate choices of the parameters $a,b,c,d$ and of the bound
states $x_{i}$ to be added. The factorization involves $m$ free
parameters defining a new tridiagonal operator $\hat{L}$, by
exchanging the two factors
\begin{equation*}
\cL=\fh(\La{a})=\cP\cQ\rightarrow \hat{\cL}=\cQ\cP=\fh(\hat{L}).
\end{equation*}
By a "rational" factorization, it is meant that each of the factors
$\mathcal{P}$ and $\mathcal{Q}$ can be written
as a polynomial in the Askey-Wilson difference operator $\La{a}$ and the
corresponding diagonal operator $\Lambda_{\gamma}$ in \eqref{1.2}, after
factoring out some rational "function" in $\Lambda_{\gamma}$.
The precise formulation is to be found in
\thref{th3.1}. The crucial \thref{th3.3} adapts the Gr\"unbaum-Yakimov
technology alluded to above, within the context of our paper.
This theorem characterizes the bispectral Darboux transformations
described in \thref{th3.1} as those for which
each factor $\cP$ and $\cQ$ is a difference operator with rational
coefficients in $q^{\gamma}$ which is moreover invariant under the involution
\begin{equation}                        \label{1.5}
I(q^{\gamma})=\frac{q^{-\gamma+1}}{abcd}\quad \mbox{and}\quad
I(E)=E^{-1}.
\end{equation}

Section 4 expresses some contiguous relations between the
Askey-Wilson functions, in the language of the discrete Darboux
transformation. This is used later in Section 6 to form in some cases the
appropriate polynomial $\mathcal{L}$ in \eqref{1.4} that needs to
be factorized, see \reref{re6.4}. Section 5 deals with the explicit
computation of the kernel of the operator $\mathcal{Q}$, from which the Darboux
factorization $\mathcal{L}=\mathcal{P}\mathcal{Q}$ is performed in Section 6.
It is quite remarkable that checking the invariance of
the operators $\cP$ and  $\cQ$ under the involution \eqref{1.5},
involves the use of classical formulas in the
theory of basic hypergeometric series, such as the
Sears' and Watson's transformation formulas.
Section 7 illustrates \thref{th6.2} on the simplest possible example,
adding one bound state off the continuous spectrum. When
$\varepsilon$ in \eqref{1.3} is put to be zero, this example
leads to orthogonal polynomials which are eigenfunctions of a
$q$-difference operator of order 4, providing a $q$-deformation
of the original Krall-Jacobi polynomials.

Section 8 of the paper requires some background in algebraic
geometry. Though it played a decisive role in our research, it is
not needed for the understanding of the earlier sections of the
paper. Its aim is to provide some
further solutions of the Askey-Wilson bispectral problem, as well
as to take a first step in putting the problem in the more general
context of the duality property satisfied by the Askey-Wilson functions,
that we established in Section 2. This requires to allow in the statement of
the problem itself for \emph{arbitrary} difference as well as
$q$-difference operators, leaving thus the context of
orthogonal polynomials.

In \cite{HI2}, we started a systematic study of
bispectral commutative rings of difference operators, for which the dual
ring is a ring of differential operators, a problem which can be thought
of as an extension of Krall's original question. The techniques we used
were adapted from the beautiful paper of G. Wilson \cite{Wi1}, who
classified all bispectral rank 1 commutative rings of
differential operators. By definition, the rank of a commutative
ring of differential or difference operators is the greatest
common divisor of the orders of all the operators in the ring.
Here again the basic philosophy is "duality": rank 1 commutative rings of
differential or difference operators enjoy some bispectral property,
when the common eigenfunction of the operators in these rings possesses a
symmetry which allows for an exchange of the geometric and the spectral
variable. Ultimately, this property reflects itself in the
geometry of the moduli space of algebraic curves which describes the
spectrum (in the sense of algebraic geometry) of these rings. In the context
of \cite{Wi1} and \cite{HI2}, the spectrum of the corresponding rings
must be an affine irreducible \emph{rational} curve with only
\emph{cusp-like singularities}.

Applying the techniques of \cite{HI2}, in \prref{pr8.8}, we exhibit a class
of \emph{rational} curves with \emph{double} points which are the spectrum of
rank 1 bispectral commutative rings of difference operators, with
dual ring a ring of $q$-difference operators. If we further impose to
these rings to contain a tridiagonal operator, we are led to
the special family of rational curves with double points given by
\begin{equation}                        \label{1.6}
v^{2}=(u^{2}-4)\prod_{i\in
J}\big(u-(q^{k_{i}/2}+q^{-k_{i}/2})\big)^{2} \prod_{i\in
K}\big(u+(q^{k_{i}/2}+q^{-k_{i}/2})\big)^{2},
\end{equation}
with $k_{i}, 1\leq i\leq g$, an arbitrary collection of positive
integers, which have been partitioned in two arbitrary sets $J$ and
$K$. This is the content of Theorem 8.1, which is the
central result of Section 8. The corresponding tridiagonal operators depend on
$g$ arbitrary parameters. They can be obtained
by iteration of the Darboux transformation, starting from the operator
\begin{equation}                        \label{1.7}
E+E^{-1},
\end{equation}
and adding bound states at the points $q^{k_{i}/2}+q^{-k_{i}/2}$
for $i\in J$, and $-(q^{k_{i}/2}+q^{-k_{i}/2})$ for $i\in K$.
It is well known that, when appropriately parametrized, the tridiagonal
operators obtained in this way provide solitonic solutions to the
doubly infinite Toda lattice hierarchy, see \cite{T}. Thus, we call these
solutions \emph{Askey-Wilson type solitons}. They are the analogues of
the rational solutions of the Korteweg-de Vries hierarchy, which were found
by Duistermaat and Gr\"unbaum \cite{DG} to form one of the two families
of solutions to the purely differential version of
the bispectral problem. When $q\to 1$, the $q$-difference bispectral operator
(appropriately scaled) in \eqref{1.2} becomes a differential
operator. All the added bound states accumulate at the boundary $\pm 2$
of the continuous spectrum $[-2,2]$ of the operator \eqref{1.7}. The
corresponding curves \eqref{1.6} acquire then only cusp-like singularities.
We refer the reader to \cite{HI3} for a complete discussion of this
limiting case.

In conclusion, the solutions of the Askey-Wilson bispectral problem
described in \thref{8.1}, as well as those described in \thref{th6.2} and
\coref{co7.6}
parallel the so-called rank 1 and rank 2 solutions found by Duistermaat and
Gr\"unbaum \cite{DG}, for the purely differential version of the
problem. Whether we have found all solutions of the Askey-Wilson bispectral
problem can only be pure speculation at this point. We feel that
finding an appropriate "moduli space" of basic hypergeometric
functions on which our result would follow from duality, looks to
be an extremely interesting and challenging problem.
Let us only mention that, even within the context of the Askey-Wilson
polynomials, the deeper understanding of duality stems from affine Hecke
algebraic considerations (see \cite{NS}), which were pioneered by Cherednik
\cite{C1}, \cite{C2} in his proof of some conjectures about the
Macdonald polynomials.

\section{Askey-Wilson functions}\label{se2}

Throughout the paper we use the standard notations for basic
hypergeometric series, following the book \cite{GR}
by Gasper and Rahman. In particular, we write
\begin{equation*}
(a_1,a_2,\dots,a_r;q)_k=\prod_{i=1}^r(a_i;q)_k,
\end{equation*}
with
\begin{equation*}
(a;q)_k=\frac{(a;q)_{\infty}}{(aq^k;q)_{\infty}}\quad\text{and}\quad
(a;q)_{\infty}=\prod_{i=0}^{\infty}(1-aq^i),
\end{equation*}
for products of $q$-shifted factorials, where $0<q<1$.
The series expansion
\begin{equation*}
{}_{r}\phi_s\left[\begin{matrix} a_1, a_2,\dots, a_{r}\\
b_1,b_2,\dots,b_s \end{matrix}\,; q,z\right]=
\sum_{k=0}^\infty\frac{\left(a_1,a_2,\dots,a_{r};q\right)_k}
{\left(q,b_1,\dots,b_s;q\right)_k}
\left[(-1)^kq^{\binom{k}{2}}\right]^{1+s-r}z^k
\end{equation*}
defines the ${}_{r}\phi_s$ basic hypergeometric series.

We denote by $E$ and $E^{-1}$, respectively,  the customary forward and
backward shift operators, acting on a function $h_{\gamma}=h(\gamma)$ by
\begin{equation*}
E h_{\gamma}=h_{\gamma +1},\qquad E^{-1} h_{\gamma}=h_{\gamma -1}.
\end{equation*}
The Askey-Wilson difference operator $\Lg$ is the second-order
difference operator
\begin{equation}                                                 \label{2.1}
\Lg=A_{\gamma} E+ B_{\gamma} \Id + C_{\gamma} E^{-1},
\end{equation}
where $\Id$ is the identity operator and the coefficients
$A_{\gamma}=\Ag{\gamma}$, $B_{\gamma}=\Bg$, $C_{\gamma}=\Cg{\gamma}$ are
given by
\begin{subequations}                        \label{2.2}
\begin{align}
A_{\gamma}&=\frac{(1-abq^{\gamma})(1-acq^{\gamma})(1-adq^{\gamma})
(1-abcdq^{\gamma-1})}{a(1-abcdq^{2\gamma-1})(1-abcdq^{2\gamma})},\label{2.2a}\\
C_{\gamma}&=\frac{a(1-q^{\gamma})(1-bcq^{\gamma-1})(1-bdq^{\gamma-1})
(1-cdq^{\gamma-1})}{(1-abcdq^{2\gamma-2})(1-abcdq^{2\gamma-1})}, \label{2.2b}\\
B_{\gamma}&=a+a^{-1}-(A_{\gamma}+C_{\gamma}).           \label{2.2c}
\end{align}
\end{subequations}

Two linearly independent solutions of the functional equation
\begin{equation}                        \label{2.3}
\Lg h_{\gamma}(z)= (z+1/z)h_{\gamma}(z)
\end{equation}
in terms of very-well-poised ${}_8\phi_7$ series, were constructed by Ismail and
Rahman in \cite{IR}, \cite{R2}. They are given by
\begin{align}
\rga&=\frac{(abq^{\gamma},acq^{\gamma},adq^{\gamma},
bcdq^{\gamma}/z;q)_{\infty}}{(bcq^{\gamma},bdq^{\gamma},cdq^{\gamma},
azq^{\gamma};q)_{\infty}}\left(\frac{a}{z}\right)^{\gamma}\times\nonumber \\
&\qquad\Wes{bcd/qz}{b/z}{c/z}{d/z}{abcdq^{\gamma-1}}{q^{-\gamma}}{qz/a},
                                \label{2.4}\\
\intertext{and}
\sga&=\frac{(abcdq^{2\gamma},bzq^{\gamma+1},czq^{\gamma+1},
dzq^{\gamma+1},bcdzq^{\gamma};q)_{\infty}}{(bcq^{\gamma},bdq^{\gamma},
cdq^{\gamma},q^{\gamma+1},bcdzq^{2\gamma+1};q)_{\infty}}(az)^{\gamma}
                            \times\nonumber\\
&\qquad\Wes{bcdzq^{2\gamma}}{bcq^{\gamma}}{bdq^{\gamma}}{cdq^{\gamma}}
{q^{\gamma+1}}{zq/a}{az},                   \label{2.5}
\end{align}
where $\Wes{a}{a_1}{a_2}{a_3}{a_4}{a_5}{z}$ stands for the very-well-poised
hypergeometric series
\begin{equation*}
{}_8\phi_7\left[\begin{matrix}a,\; qa^{1/2},\; -qa^{1/2},\; a_1,\; a_2,\;
a_3,\; a_4,\; a_5\\
a^{1/2}, -a^{1/2}, qa/a_1, qa/a_2, qa/a_3, qa/a_4, qa/a_5\end{matrix}
\,;q,z\right].
\end{equation*}
Let us denote by $D_z$ and $D_z^{-1}$, respectively, the forward and backward
$q$-shift operators, acting on a function $h(z)$ by
\begin{equation*}
D_z h(z)=h(qz)\quad\text{and}\quad D_z^{-1}h(z)=h(z/q).
\end{equation*}
Finally, let $\Bz$ be the Askey-Wilson second-order $q$-difference
operator
\begin{equation}                        \label{2.6}
\Bz = A(z)D_z-\left[A(z)+A(1/z)\right]\Id +A(1/z)D_z^{-1},
\end{equation}
where
\begin{equation}                        \label{2.7}
A(z)=\Az=\frac{(1-az)(1-bz)(1-cz)(1-dz)}{(1-z^2)(1-qz^2)}.
\end{equation}
With these notations we have
\begin{theorem}\label{th2.1}                    
The functions
\begin{subequations}                        \label{2.8}
\begin{align}
R_{\gamma}(a,b,c,d;z)&=\frac{(az;q)_{\infty}}{(bcd/z;q)_{\infty}}\rga
                                \label{2.8a}\\
\intertext{and}
S_{\gamma}(a,b,c,d;z)&=\frac{z^{-1+\log(abcd)/\log(q)}(az;q)_{\infty}}
{(zq/a,zq/b,zq/c,zq/d;q)_{\infty}}\sga              \label{2.8b}
\end{align}
\end{subequations}
solve the bispectral problem
\begin{align}
\Lg\Psi(\gamma,z)&=(z+1/z)\Psi(\gamma,z),           \label{2.9}\\
\Bz\Psi(\gamma,z)&=\Lambda_{\gamma}\Psi(\gamma,z),      \label{2.10}
\end{align}
where
\begin{equation}                        \label{2.11}
\Lambda_{\gamma}=q^{-\gamma}(1-q^{\gamma})(1-abcdq^{\gamma-1}).
\end{equation}
\end{theorem}
Note that $R_{\gamma}(a,b,c,d;z)$ and $S_{\gamma}(a,b,c,d;z)$ differ from
$\rga$ and $\sga$, respectively, by a factor independent of $\gamma$, hence
\eqref{2.9} is automatically satisfied. The proof of \eqref{2.10} is based
on an interesting bispectral involution. As we shall see, when appropriately
normalized, $\rga$ and $\sga$ satisfy an important duality relation. All
this is the content of the next lemma.

Let us introduce ``dual'' parameters and variables via the formulas
\begin{subequations}                        \label{2.12}
\begin{gather}
\at=\sqrt{\frac{q^3}{abcd}},    \quad \bt=\frac{q^2}{\at ab},\quad
\ct=\frac{q^2}{\at ac},     \quad\dt=\frac{q^2}{\at ad},    \label{2.12a}\\
\zt=\frac{q^{\gamma+1}}{\at},\quad\gt=\frac{\log(az)}{\log(q)}-1,\label{2.12b}
\end{gather}
\end{subequations}
and define
\begin{align}
{\bar{r}}(a,b,c,d;z)&=bcdz^{2-\log(abcd)/\log(q)}
\frac{(bc,bd,cd,qz/a,qz/b,qz/c,qz/d;q)_{\infty}}{(bz,cz,dz,bcd/z;q)_{\infty}},
                                \label{2.13}\\
{\bar{s}}(a,b,c,d;z)&=\frac{bcdz}{(bz,cz,dz;q)_{\infty}}.   \label{2.14}
\end{align}
Denote by $\bar{R}$ and $\bar{S}$ the functions
\begin{align}
{\bar{R}}_{\gamma}(a,b,c,d;z)&={\bar{r}}(a,b,c,d;z)\rga,    \label{2.15}\\
{\bar{S}}_{\gamma}(a,b,c,d;z)&={\bar{s}}(a,b,c,d;z)\sga.    \label{2.16}
\end{align}

\begin{lemma}\label{le2.2}                  
The functions $\bar{R}$ and $\bar{S}$ defined above satisfy the duality
relations
\begin{align}
{\bar{R}}_{\gamma}(a,b,c,d;z)&={\bar{R}}_{\gt}(\at,\bt,\ct,\dt;\zt),
                                \label{2.17}\\
\intertext{and}
{\bar{S}}_{\gamma}(a,b,c,d;z)&={\bar{S}}_{\gt}(\at,\bt,\ct,\dt;\zt).
                                \label{2.18}
\end{align}
\end{lemma}
\begin{proof}[Proof of \leref{le2.2}] Both \eqref{2.17} and \eqref{2.18} can
be proved by using transformation formula \cite[(III.24), p.~243]{GR}.
Indeed, replacing $a$, $b$, $c$, $d$, $e$, $f$ by $bcd/qz$,
$abcdq^{\gamma-1}$, $c/z$, $d/z$, $b/z$, $q^{-\gamma}$ we have
\begin{align}
&\Wes{bcd/qz}{b/z}{c/z}{d/z}{abcdq^{\gamma-1}}{q^{-\gamma}}{qz/a}\nonumber\\
&\quad =
\frac{(bcd/z,abcdq^{\gamma-1},qb/a,qc/a,qd/a,zq^{-\gamma+1}/a;q)_{\infty}}{(bc,
bd,cd,bcdq^{\gamma}/z,q^{-\gamma+2}/a^2,qz/a;q)_{\infty}}\times \label{2.19}\\
&\Wes{q^{-\gamma+1}/a^2}{q^{-\gamma+1}/ab}{q^{-\gamma+1}/ac}
{q^{-\gamma+1}/ad}{qz/a}{q/az}{abcdq^{\gamma-1}}.       \nonumber
\end{align}
The ${}_8W_7$ series in the right-hand side of \eqref{2.19} can be rewritten
in terms of the dual parameters (see \eqref{2.12}) as
\begin{equation*}
\Wes{\bt\ct\dt/q\zt}{\bt/\zt}{\ct/\zt}{\dt/\zt}{\at\bt\ct\dt q^{\gt-1}}
{q^{-\gt}}{q\zt/\at},
\end{equation*}
which is exactly the ${}_8W_7$ factor in $r_{\gt}(\at,\bt,\ct,\dt;\zt)$.
Now taking into account \eqref{2.4}, \eqref{2.12}, \eqref{2.13} and
\eqref{2.19}, one can easily check that
\begin{equation*}
\frac{\rga}{r_{\gt}(\at,\bt,\ct,\dt;\zt)}=
\frac{{\bar{r}}(\at,\bt,\ct,\dt;\zt)}{{\bar{r}}(a,b,c,d,z)},
\end{equation*}
which combined with \eqref{2.15} gives \eqref{2.17}. The proof of
\eqref{2.18} is similar. First we apply \cite[(III.24), p.~243]{GR} with
$a$, $b$, $c$, $d$, $e$, $f$ replaced by $bcdzq^{2\gamma}$, $q^{\gamma+1}$,
$bcq^{\gamma}$, $bdq^{\gamma}$, $cdq^{\gamma}$, $zq/a$ to obtain
\begin{align}
&\Wes{bcdzq^{2\gamma}}{bcq^{\gamma}}{bdq^{\gamma}}{cdq^{\gamma}}
{q^{\gamma+1}}{zq/a}{az}                    \nonumber\\
&\quad =
\frac{(bcdzq^{2\gamma+1},q^{\gamma+1},abczq^{\gamma},
abdzq^{\gamma},acdzq^{\gamma},qz^2;q)_{\infty}}
{(bzq^{\gamma+1},czq^{\gamma+1},dzq^{\gamma+1},abcdq^{2\gamma},
abcdz^2q^{\gamma},az;q)_{\infty}}\times             \label{2.20}\\
&\quad\Wes{abcdz^2q^{\gamma-1}}{az}{bz}{cz}{dz}{abcdq^{\gamma-1}}
{q^{\gamma+1}}.                         \nonumber
\end{align}
Again the hypergeometric series in the right-hand side can be
written in the dual parameters as
\begin{equation*}
\Wes{\bt\ct\dt\zt q^{2\gt}}{\bt\ct q^{\gt}}{\bt\dt q^{\gt}}{\ct\dt q^{\gt}}
{q^{\gt+1}}{\zt q/\at}{\at\zt},
\end{equation*}
and the proof follows by combining \eqref{2.5}, \eqref{2.20},
\eqref{2.12}, \eqref{2.14} and \eqref{2.16}.
\end{proof}

\begin{proof}[Proof of \thref{th2.1}]
Let $F=\Fg{\gamma}$ be ${\bar{R}}_{\gamma}(a,b,c,d)$
or ${\bar{S}}_{\gamma}(a,b,c,d)$. From \leref{le2.2}, it follows that $F$
satisfies also a difference equation in $\gt$, i.e. we have
\begin{equation}                        \label{2.21}
L_{\at,\bt,\ct,\dt}(\gt,E_{\gt})\Fgt{\gt}=(\zt+1/\zt)\Fgt{\gt}.
\end{equation}
On the other hand, from \eqref{2.12b}, it is clear that
\begin{equation}                        \label{2.22}
\Fgt{\gt\pm 1}=D_z^{\pm 1}\Fg{\gamma}.
\end{equation}
Combining \eqref{2.1}, \eqref{2.2c}, \eqref{2.21} and \eqref{2.22} we see
that $\Fz{z}$ satisfies the following $q$-difference equation in $z$
\begin{equation}                        \label{2.23}
\begin{split}
&\at(\zt+1/\zt-\at-1/\at)\Fz{z}\\
&\quad= {\bar{A}}(z)\left(\Fz{zq}-\Fz{z}\right)\\
&\qquad\quad{\bar{C}}(z)\left(\Fz{z/q}-\Fz{z}\right),
\end{split}
\end{equation}
where
${\bar{A}}(z)=\at A_{\gt}(\at,\bt,\ct,\dt)$ and
${\bar{C}}(z)=\at C_{\gt}(\at,\bt,\ct,\dt)$.
Using \eqref{2.2a}, \eqref{2.2b}, and \eqref{2.12}, one can express the
coefficients of the $q$-difference operator and the eigenvalue in formula
\eqref{2.23} in terms of the parameters $a$, $b$, $c$, $d$ and of the
variables $z$ and $\gamma$ as
\begin{subequations}                        \label{2.24}
\begin{align}
{\bar{A}}(z)&=
\frac{(1-qz/a)(1-qz/b)(1-qz/c)(1-qz/d)}{(1-z^2)(1-qz^2)},   \label{2.24a}\\
{\bar{C}}(z)&=\frac{q^3}{abcd}
\frac{(1-az/q)(1-bz/q)(1-cz/q)(1-dz/q)}{(1-z^2)(1-z^2/q)},  \label{2.24b}
\end{align}
\end{subequations}
and
\begin{equation}                        \label{2.25}
\at (\zt+1/\zt-\at-1/\at)
=q^{\gamma+1}+\frac{q^{2-\gamma}}{abcd}-\frac{q^3}{abcd}-1.
\end{equation}
If we define $\Phi(z)$ to be the function
\begin{equation}                        \label{2.26}
\Phi(z):=z^{-2+\log(abcd)/\log(q)}\frac{(az,bz,cz,dz;q)_{\infty}}
{(zq/a,zq/b,zq/c,zq/d;q)_{\infty}},
\end{equation}
then from \eqref{2.23}, \eqref{2.24} and \eqref{2.25} it follows that
$G(z):=\Phi(z)\Fz{z}$ satisfies the $q$-difference equation
\begin{equation}                        \label{2.27}
\left(q^{\gamma+1}+\frac{q^{2-\gamma}}{abcd}-\frac{q^3}{abcd}-1\right)G(z)
={\hat{A}}(z)G(zq)-{\hat{B}}(z)G(z)+{\hat{C}}(z)G(z/q),
\end{equation}
where
\begin{subequations}                        \label{2.28}
\begin{align}
{\hat{A}}(z)&=\frac{\Phi(z)}{\Phi(zq)}{\bar{A}}(z)=\frac{q^2}{abcd}A(z),
                                \label{2.28a}\\
{\hat{C}}(z)&=\frac{\Phi(z)}{\Phi(z/q)}{\bar{C}}(z)=\frac{q^2}{abcd}A(1/z),
                                \label{2.28b}\\
{\hat{B}}(z)&={\bar{A}}(z)+{\bar{C}}(z),            \label{2.28c}
\end{align}
\end{subequations}
with $A(z)$ defined in \eqref{2.7}. A direct computation now shows that
\begin{equation}                        \label{2.29}
{\hat{B}}(z)=\frac{q^2}{abcd}
\left(A(z)+A(1/z)+q-1-\frac{abcd}{q}+\frac{abcd}{q^2}\right),
\end{equation}
which combined with \eqref{2.27} and \eqref{2.28} proves that $G(z)$ satisfies
the Askey-Wilson $q$-difference equation \eqref{2.10}. To finish the proof,
it is enough to notice that
\begin{align}
\Phi(z){\bar{R}}_{\gamma}(a,b,c,d)&=bcd(bc,bd,cd;q)_{\infty}\,
R_{\gamma}(a,b,c,d;z),                      \label{2.30}\\
\intertext{and}
\Phi(z){\bar{S}}_{\gamma}(a,b,c,d)&=bcd\,S_{\gamma}(a,b,c,d;z), \label{2.31}
\end{align}
i.e. $\Phi(z){\bar{R}}$ and $\Phi(z){\bar{S}}$ differ from $R$ and $S$,
respectively, by unessential constant factors.
\end{proof}

\section{Bispectral Darboux transformations}\label{se3}

We denote by
\begin{equation}                        \label{3.1}
\cB=\left\langle\La{a},\Lambda_{\gamma}\right\rangle,
\end{equation}
the algebra of difference operators generated by $\La{a}$ and
$\Lambda_{\gamma}$, defined by \eqref{2.1} and \eqref{2.11}, respectively.
Similarly
\begin{equation}                        \label{3.2}
\cB'=\left\langle z+z^{-1},\Bg\right\rangle,
\end{equation}
will denote the algebra of $q$-difference operators generated by the operator
of multiplication by $z+z^{-1}$ and the operator $\Bg$ defined in \eqref{2.6}.
Formulas \eqref{2.9} and \eqref{2.10} serve to define an anti-isomorphism
\begin{equation}                        \label{3.3}
\fb:\cB\rightarrow\cB'
\end{equation}
between these two algebras, i.e. it is given on the generators by
\begin{equation}                        \label{3.4}
\fb(\La{a})=z+z^{-1}\quad\text{and}\quad \fb(\Lambda_{\gamma})=\Bg.
\end{equation}
Note that the bispectral property \eqref{2.9}, \eqref{2.10} is
equivalent to the identity
\begin{equation} \label{3.5}
X \Psi(\gamma,z)=\fb(X)\Psi(\gamma,z),\quad \forall X\in\cB.
\end{equation}
We shall also need the commutative subalgebras (the algebras of ``functions'')
of $\cB$ and $\cB'$ defined by
\begin{equation}                        \label{3.6}
\cK=\left\langle\Lambda_{\gamma}\right\rangle\quad\text{and}\quad
\cK'=\left\langle z+z^{-1}\right\rangle.
\end{equation}

The next theorem summarizes the technology of bispectral Darboux
transformations, initiated in \cite{BHY1}, \cite{BHY2} and \cite{KR}. It was
adapted and applied to the case of difference operators in \cite{GHH}.
For the convenience of the reader, we include the short proof of
this result.
\begin{theorem}\label{th3.1}                           
Let $\cL$ be a constant coefficient polynomial in $\La{a}$, which
factorizes rationally as
\begin{equation}                        \label{3.7}
\cL=\cP\cQ,
\end{equation}
in such a way that
\begin{equation}                        \label{3.8}
\cP=U\Gamma^{-1},\qquad \cQ=\Theta^{-1}V,
\end{equation}
with $U,V\in\cB$, and $\Theta, \Gamma\in\cK$. Then the Darboux transform
of $\cL$ given by
\begin{equation}                        \label{3.9}
\hat{\cL}=\cQ\cP
\end{equation}
is a bispectral operator. More precisely, defining $f=\fb(\cL)\in\cK'$ and
$\hat{\Psi}=\cQ\Psi$, with $\Psi$ satisfying \eqref{2.9} and \eqref{2.10},
we have
\begin{align}
\hat{\cL}\hat{\Psi}&=f\hat{\Psi},             \label{3.10}\\
\hat{B}\hat{\Psi}&=\Theta\Gamma\hat{\Psi},          \label{3.11}
\end{align}
with
\begin{equation}                        \label{3.12}
\hat{B}=\fb(V)\fb(U)f^{-1}.
\end{equation}
\end{theorem}
\begin{proof}
From the definitions above and equation \eqref{3.5}, we obtain
\begin{equation*}
\hat{\cL}\hat{\Psi}=\cQ\cP\cQ\Psi=\cQ\cL\Psi=\cQ\fb({\cL})\Psi=
f\cQ\Psi=f\hat{\Psi},
\end{equation*}
which establishes \eqref{3.10}. Using \eqref{3.5} and
\eqref{3.8}, we can write $\hat{\Psi}$ as
\begin{equation}                        \label{3.13}
\hat{\Psi}=\cQ\Psi=\Theta^{-1} V\Psi=\Theta^{-1}\fb(V)\Psi.
\end{equation}
From \eqref{3.7} and \eqref{3.8}, we deduce that
\begin{equation*}
\Theta\Gamma =V \cL^{-1} U.
\end{equation*}
Applying the anti-isomorphism $\fb$ to this equation, we obtain
\begin{equation}                        \label{3.14}
\fb(\Theta\Gamma)=\fb(U)f^{-1}\fb(V).
\end{equation}
Finally, using \eqref{3.5}, \eqref{3.12}, \eqref{3.13} and \eqref{3.14} we get
\begin{align*}
\Theta\Gamma\hat{\Psi}&=\Theta\Gamma\Theta^{-1}\fb(V)\Psi
=\Theta^{-1}\fb(V)\Theta\Gamma\Psi
=\Theta^{-1}\fb(V)\fb(\Theta\Gamma)\Psi\\
&=\Theta^{-1}\fb(V)\fb(U)f^{-1}\fb(V)\Psi
=\fb(V)\fb(U)f^{-1}\Theta^{-1}\fb(V)\Psi=
\hat{B}\hat{\Psi},
\end{align*}
which gives \eqref{3.11} and completes the proof.
\end{proof}

Despite the apparent simplicity of \thref{th3.1}, it is a priori
very complicated to recognize that an operator admits a rational
Darboux factorization as defined by \eqref{3.7}, \eqref{3.8}. The
rest of the section is concerned with making \thref{th3.1} more
effective. Some of the ideas developed below were first introduced in
\cite{GY}, in the context of the Jacobi polynomials.

Let us denote by $\RE$ the algebra of difference operators of the form
\begin{equation*}
T=\sum_{j=m_1}^{m_2}h_j(q^{\gamma})E^j,
\end{equation*}
with coefficients $h_j(q^{\gamma})$ rational functions in $q^{\gamma}$.
The ordered pair $[m_1,m_2]$ is called the support of $T$.
Consider the involution\footnote{By involution we mean an automorphism $I$,
such that $I\circ I=\Id$.} $I$ on $\RE$, defined by
\begin{equation}                        \label{3.15}
I(q^{\gamma})=\frac{q^{-\gamma+1}}{abcd}\quad\text{and}\quad I(E)=E^{-1}.
\end{equation}
\begin{lemma} \label{le3.2}                              
A Laurent polynomial in $q^{\gamma}$
\begin{equation}\label{3.16}
p(q^{\gamma})=\sum_{k=-m}^{n}c_{k}q^{k\gamma},\quad
m,n\in\mathbb{N},
\end{equation}
is a polynomial in $\Lambda_{\gamma}$ (i.e. it
belongs to $\cK$ as defined in \eqref{3.6}) if and only if it is
$I$-invariant.
\end{lemma}
\begin{proof} A straightforward computation from \eqref{2.11}
gives
\begin{equation*}
I(\Lambda_{\gamma})=\Lambda_{\gamma},
\end{equation*}
hence every polynomial in $\Lambda_{\gamma}$ is $I$-invariant.
Conversely, if $p(q^{\gamma})$ in \eqref{3.16} is $I$-invariant,
one checks that $m=n$ and
\begin{equation}\label{3.17}
p(q^{\gamma})=\sum_{k=0}^{n}c_{k}\Big(q^{k\gamma}+\frac{q^{k(1-\gamma)}}{(abcd)^{k}}\Big).
\end{equation}
Observe from \eqref{2.11} that
\begin{equation*}
q^{\gamma}+\frac{q^{1-\gamma}}{abcd}=\frac{q}{abcd}\Big(\Lambda_{\gamma}+1+\frac{abcd}{q}\Big).
\end{equation*}
Since
\begin{multline*}
q^{k\gamma}+\frac{q^{k(1-\gamma)}}{(abcd)^{k}}=\Big(q^{\gamma}+\frac{q^{1-\gamma}}{abcd}\Big)
\Big(q^{(k-1)\gamma}+\frac{q^{(k-1)(1-\gamma)}}{(abcd)^{k-1}}\Big)\\-\frac{q}{abcd}
\Big(q^{(k-2)\gamma}+\frac{q^{(k-2)(1-\gamma)}}{(abcd)^{k-2}}\Big),
\end{multline*}
by induction, it follows that $p(q^{\gamma})$ in \eqref{3.17}
is a polynomial in $\Lambda_{\gamma}$, i.e. it belongs to $\cK$.
\end{proof}

An easy computation, using formulas \eqref{2.2a} and \eqref{2.2b},  gives
\begin{equation}                        \label{3.18}
I(\Aa{a})=\Ca{a},
\end{equation}
which shows that the operator $\La{a}$ in
\eqref{2.1} is $I$-invariant. Let us denote by $\dI$ the subalgebra of $\RE$,
consisting of $I$-invariant operators, i.e.
\begin{equation*}
\dI=\left\{T\in\RE|\,I(T)=T\right\}.
\end{equation*}
Since $\La{a}\in\dI$ and $\Lambda_{\gamma}\in\dI$, it follows that
\begin{equation}                        \label{3.19}
\cK\subset\cB\subset\dI.
\end{equation}
Let $\dI_m$ be the subset of $\dI$, consisting of $I$-invariant operators
with support $[-m,m]$. The next theorem characterizes the bispectral Darboux
transformations as those for which the two factors $\cP$ and $\cQ$
in \eqref{3.7} are $I$-invariant operators, with rational
coefficients.

\begin{theorem}\label{th3.3}                    
The following conditions on an operator $T\in\RE$ are equivalent:
\begin{itemize}
\item[(i)]   The operator $T$ is $I$-invariant, i.e. $T\in\dI$;
\item[(ii)]  $T$ has the form $\Theta^{-1}V$, for some operator $V\in\cB$ and
        some function $\Theta\in\cK$;
\item[(iii)] $T$ has the form $U\Gamma^{-1}$, for some operator $U\in\cB$ and
        some function $\Gamma\in\cK$.
\end{itemize}
\end{theorem}
\begin{proof} From \eqref{3.19} it is obvious that (ii) implies (i). Below, we shall prove
that (i) implies (ii), i.e. if $T\in\dI$, then $T={\Theta}^{-1}V$, for some
$\Theta\in\cK$ and $V\in\cB$. The proof of the fact that (i) is equivalent to
(iii) is similar.

It is clear that $T$ must have support $[-m,m]$, for some
$m\in\N$. If $m=0$, we can always assume that $T$ is the quotient
of two $I$-invariant Laurent polynomials in $q^{\gamma}$. Hence,
by \leref{le3.2}, $T$ is a rational function in $\Lambda_{\gamma}$, so
we can take $V$ and $\Theta$ in (ii) to be the numerator and the denominator,
respectively. If $m\geq 1$, we shall show that
\begin{equation}                        \label{3.20}
T=\Theta_{m}^{-1}V_m \quad \mbox{mod}\;\dI_{m-1},
\end{equation}
with $\Theta_{m}\in\cK$, $V_{m}\in\cB$, from which the proof follows by
induction.

Any operator in $\dI_m$ may be written
\begin{equation}                        \label{3.21}
T=(\mbox{Id}+I)\frac{\fc(q^{\gamma})}{\fd(q^{\gamma})}E^m\quad \mbox{mod}\;\dI_{m-1},
\end{equation}
with $\fc(q^{\gamma})$ and $\fd(q^{\gamma})$ polynomials in $q^{\gamma}$. In
particular, since $\La{a}$ is $I$-invariant, we have
\begin{equation}                         \label{3.22}
\La{a}^{m}=(\mbox{Id}+I)\frac{\ff(q^{\gamma})}{\fg(q^{\gamma})}E^{m}\quad
\mbox{mod}\;\dI_{m-1},
\end{equation}
for some polynomials $\ff(q^{\gamma})$ and $\fg(q^{\gamma})$.

We denote by $\lambda_{x}$ and $\rho_{x}$ the operators of left and right
multiplication by $x$, respectively, on $\RE$, i.e. for $\cT\in\RE$ we have
\begin{equation*}
\lambda_{x}(\cT)=x\cT\quad\text{and}\quad\rho_{x}(\cT)=\cT x.
\end{equation*}
Eliminating $q^{-\gamma}$ (resp. $q^{\gamma}$) from the equations
\begin{align*}
\Lambda_{\gamma}E^{m}&=(abcdq^{\gamma-1}+q^{-\gamma}-1-abcdq^{-1})E^{m},\\
E^{m}\Lambda_{\gamma}&=(q^{m}abcdq^{\gamma-1}+q^{-m}q^{-\gamma}-1-abcdq^{-1})E^{m},
\end{align*}
one obtains
\begin{equation} \label{3.23}
\delta_m^{+} (E^{m})=q^{\gamma}E^{m}\quad \mbox{and}\quad
\delta_m^{-}(E^{m})=q^{-\gamma}E^{m},
\end{equation}
with
\begin{subequations}                        \label{3.24}
\begin{align}
\delta_m^{+} &=\frac{q}{abcd(1-q^{2m})}(\lL-q^m\rL)+
\left(1+\frac{q}{abcd}\right)\frac{1}{1+q^m}\Id,        \label{3.24a}\\
\intertext{and}
\delta_m^{-}&=\frac{1}{1-q^{-2m}}(\lL-q^{-m}\rL)+
\frac{1}{1+q^{-m}}\left(\frac{abcd}{q}+1\right)\Id \nonumber\\
&=\frac{abcd}{q}\delta_{-m}^{+}.    \label{3.24b}
\end{align}
\end{subequations}

Notice that $\delta_m^{+}$ and $\delta_m^{-}$ preserve the
subalgebra $\dI$. Hence,
from \eqref{3.22}, using the properties \eqref{3.23} of
$\delta_m^{+}$ and $\delta_m^{-}$, one
deduces that
\begin{equation*}
p(\delta_{m}^{+},\delta_m^{-})\La{a}^{m}=(\mbox{Id}+I)p(q^{\gamma},q^{-\gamma})
\frac{\ff(q^{\gamma})}{\fg(q^{\gamma})}E^{m}\quad
\mbox{mod}\;\dI_{m-1},
\end{equation*}
for any polynomial $p(q^{\gamma},q^{-\gamma})\in
\mathbb{C}[q^{\gamma},q^{-\gamma}]$. In particular, if we choose
\begin{equation*}
p(q^{\gamma},q^{-\gamma})=\fc(q^{\gamma})\;\fg(q^{\gamma})\;I\big(\fd(q^{\gamma})\;\ff(q^{\gamma})\big),
\end{equation*}
we get
\begin{equation*}
p(\delta_m^{+},\delta_m^{-})\La{a}^{m}=\fd(q^{\gamma})\;\ff(q^{\gamma})
\;I\big(\fd(q^{\gamma})\;\ff(q^{\gamma}))\;(\mbox{Id}+I)
\frac{\fc(q^{\gamma})}{\fd(q^{\gamma})}E^m \quad \mbox{mod}\;\dI_{m-1}.
\end{equation*}
Define
\begin{equation*}
\Theta_m=\fd(q^{\gamma})\;\ff(q^{\gamma})
\;I\big(\fd(q^{\gamma})\;\ff(q^{\gamma})).
\end{equation*}
Since $\Theta_m$ is $I$-invariant, by \leref{le3.2}, $\Theta_m
\in\cK$. Also, it is clear from the definitions \eqref{3.24a} and
\eqref{3.24b} of $\delta_m^{+}$ and $\delta_m^{-}$ that $p(\delta_{m}^{+},\delta_m^{-})\La{a}^{m}
\in \cB$. Hence, remembering \eqref{3.21}, we get
\begin{align*}
\Theta_m^{-1}\;p(\delta_{m}^{+},\delta_m^{-})\La{a}^{m}&=(\mbox{Id}+I)\;
\frac{\fc(q^{\gamma})}{\fd(q^{\gamma})}\;E^m \quad
\mbox{mod}\;\dI_{m-1}\\
&=T\quad\mbox{mod}\;\dI_{m-1},
\end{align*}
which establishes \eqref{3.20} with
$V_{m}=p(\delta_{m}^{+},\delta_m^{-})\La{a}^{m}$. This completes the proof of \thref{th3.3}.
\end{proof}

\section{Contiguous relations}\label{se4}

From \eqref{2.5}, we observe that the function
$s_\gamma(aq,b,c,d;a)$ has a simple product form, since the
$_8W_7$ series defining it reduces to $1$. One checks easily that
\begin{equation*}
\frac{s_\gamma(aq,b,c,d;a)}{s_{\gamma-1}(aq,b,c,d;a)}=\frac{\varphi_\gamma}
{\varphi_{\gamma-1}}\frac{\Cg{\gamma}}{\Ag{\gamma}},
\end{equation*}
with
\begin{equation} \label{4.1}
\varphi_{\gamma}=\frac{1}{q^{-\gamma}-abcdq^{\gamma}}.
\end{equation}
Let us denote by $\Pa$ and $\Qa$ the difference operators
\begin{align}
\Pa&=\varphi_{\gamma}(E-\Id),           \label{4.2}\\
\Qa&=\left(\Ag{\gamma}\Id -\Cg{\gamma} E^{-1}\right)\,
\frac{1}{\varphi_{\gamma}}.               \label{4.3}
\end{align}
From \eqref{2.1} and \eqref{2.2c}, it follows that
\begin{equation}                        \label{4.4}
\La{a}-\left(a+a^{-1}\right)\Id=\Qa\Pa.
\end{equation}
Since the kernel of the operator $\Qa$ is generated by
$s_\gamma(aq,b,c,d;a)$, it is natural to expect that interchanging
the factors $\Qa$ and $\Pa$ in \eqref{4.4}, we obtain
the operator $\La{aq}$. This simply means that $\La{a}$ can be obtained as a
{\it Darboux transformation\/} from the operator $\La{aq}$. The
proposition below makes this statement precise.
\begin{proposition}\label{pr4.1}            
With the notations above, we have
\begin{equation}                        \label{4.5}
\La{aq}-\left(a+a^{-1}\right)\Id=\Pa\Qa.
\end{equation}
Moreover, up to a factor independent of $\gamma$, $\Qa$ maps $\ra{aq}$ and
$\sa{aq}$ into $\ra{a}$ and $\sa{a}$, respectively. More precisely,
the following contiguous relations
\begin{subequations}                        \label{4.6}
\begin{align}
\Qa\;\ra{aq}&=\frac{1-az}{a}\ra{a},             \label{4.6a}\\
\Qa\;\sa{aq}&=\frac{(z-a)(1-az)}{az}\sa{a},     \label{4.6b}
\end{align}
\end{subequations}
hold. In particular, from \eqref{4.6b}, it follows that the kernel of $\Qa$ is
spanned by $s_{\gamma}(aq,b,c,d;a)$.
\end{proposition}

\begin{proof}
A straightforward  computation shows that
\begin{equation}                        \label{4.7}
a+a^{-1}-\Cg{\gamma+1}-\Ag{\gamma}=aq+(aq)^{-1}-\Aa{aq}-\Ca{aq},
\end{equation}
which combined with \eqref{2.1}, \eqref{2.2}, \eqref{4.2}, and \eqref{4.3}
gives \eqref{4.5}. The proof of \eqref{4.6} can be easily extracted from
the contiguous relations found by Ismail and Rahman in \cite{IR}. Let us
denote for simplicity by $\rgaf{\gamma}{a}$ and $\sgaf{\gamma}{a}$ the
${}_8W_7$ factor on the right-hand side of formulas \eqref{2.4} and
\eqref{2.5}, respectively. Applying \cite[formula (2.3), p.~207]{IR}, with
$a$, $b$, $c$, $d$, $e$, $f$ replaced by $bcd/qz$, $abcdq^{\gamma}$,
$q^{-\gamma}$, $b/z$, $c/z$, $d/z$, we obtain
\begin{equation}                        \label{4.8}
\begin{split}
&\left(q^{-\gamma}-abcdq^{\gamma-1}\right)\left(1-\frac{1}{az}\right)
\rgaf{\gamma}{a}\\
&\qquad = \left(q^{-\gamma}-1\right)\left(1-\frac{bcdq^{\gamma-1}}{z}\right)
\rgaf{\gamma-1}{aq}\\
&\qquad\quad +\left(1-abcdq^{\gamma-1}\right)
\left(1-\frac{1}{azq^{\gamma}}\right)\rgaf{\gamma}{aq}.
\end{split}
\end{equation}
From \eqref{2.4} and \eqref{4.8}, one can deduce \eqref{4.6a}.
To get \eqref{4.6b}, we first substitute $bcdzq^{2\gamma+1}$, $z/a$,
$bcq^{\gamma}$, $bdq^{\gamma}$, $cdq^{\gamma}$, $q^{\gamma+1}$ for
$A^2$, $A/\lambda$, $A/\mu$, $A/\nu$, $A/\rho$, $A/\sigma$ in
\cite[formula (2.17), p.~210]{IR} to obtain
\begin{align}
&\sgaf{\gamma-1}{aq}
=\frac{(1-bcdzq^{2\gamma})(1-bcdzq^{2\gamma-1})}
{(1-bzq^{\gamma})(1-czq^{\gamma})(1-dzq^{\gamma})(1-bcdzq^{\gamma-1})}\times
                                \nonumber\\
&\quad\bigg[(1-z/a)(1-az)\sgaf{\gamma-1}{a}         \label{4.9}\\
&\quad\quad+\frac{z}{a}\frac{(1-abq^{\gamma})(1-acq^{\gamma})
(1-adq^{\gamma})(1-abcdq^{\gamma-1})}{(1-abcdq^{2\gamma-1})(1-abcdq^{2\gamma})}
\sgaf{\gamma}{aq}\bigg],                    \nonumber
\end{align}
which combined with \eqref{2.5} gives \eqref{4.6b}.
\end{proof}

\begin{remark}\label{re4.2}                     
We can also consider $\La{aq}$ as obtained by means of a Darboux
transformation from $\La{a}$. In this case, the analogous formulas to
\eqref{4.6a} and \eqref{4.6b} are
\begin{subequations}                        \label{4.10}
\begin{align}
\Pa\;\ra{a}&=\frac{a-z}{z}\ra{aq},              \label{4.10a}\\
\Pa\;\sa{a}&=-\sa{aq}.                      \label{4.10b}
\end{align}
\end{subequations}
One way to prove these formulas is to use \cite[formula (2.2), p.~207]{IR}
and \cite[formula (2.18), p.~210]{IR}, and to proceed as before, adapting
the proof of \prref{pr4.1}. However, they can be easily deduced from
\prref{4.1}, by applying $\Pa$ to both sides of \eqref{4.6a}-\eqref{4.6b},
and by using \eqref{4.5} and the fact that $\ra{a}$ and $\sa{a}$ satisfy
\eqref{2.9}.
\end{remark}

\section{Computing the kernel of the Darboux transform using Sears' and
Watson's transformation formulas}\label{se5}

As we saw in \thref{th3.3}, in order to apply the general bispectral
Darboux technology, one needs to find an $I$-invariant factorization of
a constant coefficient polynomial $\cL$ in $\La{a}$. In this section
we construct $I$-invariant eigenfunctions of $\La{a}$, which generate
the kernel of the operator $\cQ$, in the notations of \thref{th3.1}.

\begin{proposition}\label{pr5.1}            
Assume that
\begin{equation}                        \label{5.1}
a=dq^{\alpha}\quad \text{and}\quad z_{k}=dq^{k},\;k,\alpha\in \N,
\;0\leq k\leq \alpha-1.
\end{equation}
Then,
\begin{align}
r_{\gamma}(z_{k})&=(bcq^{-k},d^{2}q^{\alpha};q)_{k}\times   \nonumber\\
&\qquad f_{\gamma}\;\phift{d^{2}q^{k}}{bcd^{2}q^{\gamma+\alpha-1}}
{q^{-\gamma}}{q^{-k}}{bd}{cd}{d^{2}q^{\alpha}}{q},      \label{5.2}\\
s_{\gamma}(z_{k})&=(bc)^{\alpha-k-1}(d^{2}q^{k},dq^{k+1}/b,
dq^{k+1}/c;q)_{\alpha-k-1}\times                \nonumber\\
&\qquad f_{\gamma}g_{\gamma}\;\phift{q^{k-\alpha+1}}{bcq^{\gamma}}
{q^{-\gamma-\alpha+1}/d^{2}}{q^{1-\alpha-k}/d^{2}}{q^{2-\alpha}/d^{2}}
{bq^{1-\alpha}/d}{cq^{1-\alpha}/d}{q},              \label{5.3}
\end{align}
with
\begin{align}
f_{\gamma}&=\frac{q^{\gamma\alpha}}
    {(bdq^{\gamma},cdq^{\gamma};q)_{\alpha}},       \label{5.4}\\
g_{\gamma}&=q^{\gamma(\alpha-1)}d^{2\gamma}
\frac{(d^{2}q^{\gamma+\alpha},bcd^{2}q^{\gamma+\alpha-1};q)_{\infty}}
{(q^{\gamma+1},bcq^{\gamma};q)_{\infty}}.           \label{5.5}
\end{align}
\end{proposition}

If we make the further assumption that
\begin{equation}                        \label{5.6}
d^{2}=q^{l},\quad l\in\N,\; l\geq 1,
\end{equation}
$g_{\gamma}$ in \eqref{5.5} can be rewritten as
\begin{equation}                        \label{5.7}
g_{\gamma}=\frac{q^{\gamma(l+\alpha-1)}}
    {(q^{\gamma+1},bcq^{\gamma};q)_{l+\alpha-1}}.
\end{equation}

\begin{proof} Applying Watson's transformation formula
\cite[(III.18), p.~242]{GR} with $n$, $a$, $b$, $c$, $d$, $e$ replaced by
$k$, $bcq^{-k-1}$, $bq^{-k}/d$, $cq^{-k}/d$, $bcd^2q^{\gamma+\alpha-1}$,
$q^{-\gamma}$, respectively, we can rewrite $r_{\gamma}(z_k)$ as
\begin{align}
r_{\gamma}(z_k)&=q^{(\alpha-k)\gamma}\frac{(d^2q^{\gamma+\alpha},bcq^{-k};q)_k}
{(bdq^{\gamma},cdq^{\gamma};q)_{\alpha}}
\frac{(q^{-\alpha-k+1}/d^2;q)_k}
{(q^{-\gamma-\alpha-k+1}/d^2;q)_{k}}\times\nonumber\\
&\qquad\phift{d^2q^k}{bcd^2q^{\gamma+\alpha-1}}{q^{-\gamma}}{q^{-k}}{cd}{bd}
{d^2q^{\alpha}}{q}.                     \label{5.8}
\end{align}
Using formula (I.9) on page 233 in \cite{GR}, one can easily deduce that
\begin{equation*}
\frac{(q^{-\alpha-k+1}/d^2;q)_k}{(q^{-\gamma-\alpha-k+1}/d^2;q)_{k}}=
q^{\gamma k}\frac{(d^2q^\alpha;q)_k}{(d^2q^{\gamma+\alpha};q)_k},
\end{equation*}
which shows that \eqref{5.8} is equivalent to \eqref{5.2}. The proof of
\eqref{5.3} is similar. Let us first replace the parameters
$n$, $a$, $b$, $c$, $d$, $e$ in Watson's formula \cite[(III.18), p.~242]{GR} by
$\alpha-k-1$, $bcd^2q^{2\gamma+k}$, $bdq^{\gamma}$, $cdq^{\gamma}$,
$bcq^{\gamma}$, $q^{\gamma+1}$, to obtain the following formula for
$s_{\gamma}(z_k)$ in terms of ${}_4\phi_3$ series
\begin{align}
s_{\gamma}(z_k)&=q^{(\alpha+k)\gamma}d^{2\gamma}
\frac{(d^2q^{\gamma+\alpha},bcd^2q^{\gamma+\alpha-1};q)_{\infty}}
{(q^{\gamma+1},bcq^{\gamma};q)_{\infty}}
\frac{(d^2q^k;q)_{\alpha-k-1}}{(bdq^{\gamma},cdq^{\gamma};q)_{k+1}}
    \times                          \nonumber\\
&\qquad\phift{q^{k+1}}{bcq^{\gamma}}{q^{\gamma+1}}{q^{-\alpha+k+1}}
{cdq^{\gamma+k+1}}{bdq^{\gamma+k+1}}{q^{2-\alpha}/d^2}{q}.  \label{5.9}
\end{align}
Now replacing $n$, $a$, $b$, $c$, $d$, $f$  by $\alpha-k-1$,
$bcq^{\gamma}$, $q^{k+1}$, $q^{\gamma+1}$, $q^{2-\alpha}/d^2$,
$cdq^{\gamma+k+1}$, $bdq^{\gamma+k+1}$ in Sears' transformation formula
(see \cite[(III.15), p.~242]{GR}), we see that \eqref{5.9} gives
exactly \eqref{5.3}.
\end{proof}

In the case considered above (i.e. $a=dq^{\alpha}$ and $d^2=q^l$)
formula \eqref{3.15} reduces to
\begin{equation}                        \label{5.10}
I(q^{\gamma})=\frac{q^{-\gamma-\alpha+1}}{bcd^2}=
\frac{q^{-\gamma-\alpha-l+1}}{bc}.
\end{equation}

\begin{lemma}\label{le5.2}                  
If \eqref{5.6} holds, the functions $r_{\gamma}(z_k)$ and $s_{\gamma}(z_k)$,
defined by \eqref{5.2} and \eqref{5.3} are $I$-invariant rational functions
in $q^{\gamma}$, i.e. we have $I(r_{\gamma}(z_k))=r_{\gamma}(z_k)$ and
$I(s_{\gamma}(z_k))=s_{\gamma}(z_k)$.
\end{lemma}

\begin{proof}
From \eqref{5.2}, \eqref{5.3}, \eqref{5.4}, and \eqref{5.7}, it is obvious
that $r_{\gamma}(z_k)$ and $s_{\gamma}(z_k)$ are rational functions in
$q^{\gamma}$. Since $I(bcd^2q^{\gamma+\alpha-1})=q^{-\gamma}$ and
$I(bcq^{\gamma})=q^{-\gamma-\alpha+1}/d^2$, it is clear that the ${}_4\phi_3$
series in \eqref{5.2} and \eqref{5.3} are $I$-invariant. Thus, to show the
invariance of $r_{\gamma}(z_k)$ and $s_{\gamma}(z_k)$, it is enough to
show the invariance of $f_{\gamma}$ and $g_{\gamma}$. Using the definitions
\eqref{5.4} of $f_{\gamma}$ and \eqref{5.10} of $I$, one computes
\begin{equation*}
\begin{split}
I(f_{\gamma})&=\frac{\left(bcd^2q^{\gamma+\alpha-1}\right)^{-\alpha}}
{(1-1/cdq^{\gamma+\alpha-1})(1-1/cdq^{\gamma+\alpha-2})\cdots
(1-1/cdq^{\gamma})}\\
&\quad\times\frac{1}{(1-1/bdq^{\gamma+\alpha-1})
(1-1/bdq^{\gamma+\alpha-2})\cdots(1-1/bdq^{\gamma})}\\
&=\frac{q^{2\gamma\alpha+\alpha(\alpha-1)-\alpha(\gamma+\alpha-1)}}
{(bdq^{\gamma},cdq^{\gamma};q)_{\alpha}}=f_{\gamma}.
\end{split}
\end{equation*}
The proof of the $I$-invariance of $g_{\gamma}$ is similar.
\end{proof}
\section{Askey-Wilson type functions}\label{se6}

In this section we prove the main result of the paper. First we review
some basic facts about Darboux transformations.

\subsection{Iterating the Darboux transformation}
Consider a difference operator $L_0$ and a nonzero eigenfunction
$\psi_1(\gamma)$ for which $L_0\psi_1(\gamma)=x_1\psi_1(\gamma)$.
If $Q_0$ is a first-order difference operator with kernel spanned by
$\psi_1(\gamma)$, we can write $L_0$ in the form $L_0=x_1\Id+P_0Q_0$,
for some difference operator $P_0$. The operator $L_1=x_1\Id+Q_0P_0$ is
by definition a Darboux transformation of $L_0$.

The next proposition, which can be found in \cite[pp.~14-19]{Wr}
for the case of differential operators,
describes the result of $m$ successive Darboux transformations starting
from $L_0$
\begin{align}
&L_0=x_{j_1}\Id+P_0Q_0\curvearrowright L_1=x_{j_1}\Id+Q_0P_0=x_{j_2}\Id+P_1Q_1
\curvearrowright\cdots                      \nonumber\\
&\quad L_{m-1}=x_{j_{m-1}}\Id+Q_{m-2}P_{m-2}=x_{j_m}\Id+P_{m-1}Q_{m-1}
                                \label{6.1}\\
&\qquad\curvearrowright \hL=L_m=x_{j_m}\Id+Q_{m-1}P_{m-1}.  \nonumber
\end{align}

\begin{proposition} \label{pr6.1}                      
If the operator $\hL$ is obtained from $L_0$ by iteration
of the Darboux transformation \eqref{6.1} then
\begin{equation}                        \label{6.2}
\hL \cQ=\cQ L_0,
\end{equation}
where
\begin{equation}                        \label{6.3}
\cQ=Q_{m-1}Q_{m-2}\cdots Q_0.
\end{equation}

Conversely, if $\hL$ and $L_0$ are difference operators and if there exists
a difference operator $\cQ$ of order $\geq 1$ such that \eqref{6.2} holds,
then $\hL$ can be obtained from $L_0$ by a sequence of Darboux transformations.
\end{proposition}
\begin{proof} From \eqref{6.1} one can deduce that
\begin{equation}                        \label{6.4}
L_kQ_{k-1}=Q_{k-1}L_{k-1}\text{ for } k=1,2,\dots,m,
\end{equation}
which easily implies \eqref{6.2}, thus proving the first assertion.

The second part can be proved by induction on the order of the operator $\cQ$
as follows. Assume that \eqref{6.2} holds. It is obvious that $\Ker \cQ$ is
preserved by $L_0$, i.e. $L_0(\Ker \cQ)\subset\Ker \cQ$. Thus we may regard
$L_0$ restricted to $\Ker\cQ$ as a linear operator acting in a
finite-dimensional complex vector space. Let $\{\psi_j(\gamma)\}_{j=1}^m$ be
a Jordan basis of this operator, i.e.
\begin{subequations}
\begin{equation}                        \label{6.5a}
\Ker\cQ=\Span\{\psi_1(\gamma),\psi_2(\gamma),\dots,\psi_m(\gamma)\},
\end{equation}
where
\begin{equation}                        \label{6.5b}
L_0\psi_k(\gamma)=x_{j_k}\psi_k(\gamma)+\sigma_k\psi_{k-1}(\gamma),
\qquad 1\leq k \leq m,
\end{equation}
\end{subequations}
with $\sigma_1=0$ and $\sigma_k=0$ or $1$ for $2\leq k\leq m$. Let $Q_0$ be
a first-order difference operator with kernel spanned by $\psi_1(\gamma)$.
If the order of $\cQ$ is 1, we simply take $Q_0=\cQ$. Otherwise we can take
for example
\begin{equation*}
Q_0=E-\frac{\psi_1(\gamma+1)}{\psi_1(\gamma)}\Id.
\end{equation*}
Since $\psi_1(\gamma)$ is in the kernels of $L_0-x_{j_1}\Id$ and $\cQ$,
there exist difference operators $P_0$ and $\tQ$ such that
\begin{align}
L_0-x_{j_1}\Id=P_0Q_0                    \label{6.6}
\intertext{ and }
\cQ=\tQ Q_0.                             \label{6.7}
\end{align}
If we define $L_1=x_{j_1}\Id+Q_0P_0$, then from \eqref{6.6} we get that
\begin{equation}                         \label{6.8}
L_1Q_0=Q_0L_0.
\end{equation}
Using the last equation, \eqref{6.2} and \eqref{6.7} one can show that
\begin{equation*}
\hL\tQ=\tQ L_1.
\end{equation*}
Moreover from \eqref{6.5b} and \eqref{6.8} we can see that for $k=2,3,\dots,m$
we have
\begin{equation*}
L_1(Q_0\psi_k)=x_{j_k}Q_0\psi_k+\sigma'_kQ_0\psi_{k-1}, \text{ where }
\sigma'_2=0, \text{ and }\sigma'_k=\sigma_k \text{ if }k>2,
\end{equation*}
i.e. $\{Q_0\psi_j(\gamma)\}_{j=2}^m$ is a Jordan basis for the operator
$L_1$ in the vector space $\Ker\tQ$. The proof now follows by induction.
\end{proof}

The above proposition tells us that the Darboux process \eqref{6.1} is
determined by the intertwining operator $\cQ$ in \eqref{6.2}. The points
$\{x_{j_k}\}$ are simply the eigenvalues of the operator $L_0$ in the vector
space $\Ker \cQ$. We will use \prref{pr6.1} as follows. First, we will
construct a difference operator $\cQ$ whose kernel is $L_0$ invariant, with
$L_0=\La{a}$ the Askey-Wilson operator defined in \seref{se2}. This operator
$\cQ$ will correspond to a sequence of Darboux transformations \eqref{6.1}.
Indeed if we consider the operator $\cQ L_0$, then its kernel contains the
kernel of $\cQ$, so it has a factorization of the form
\begin{equation*}
\hL \cQ=\cQ L_0,
\end{equation*}
for some (unique) difference operator $\hL$. By the above considerations,
$\hL$ is obtained from $L_0$ by a Darboux chain \eqref{6.1}.

The next important fact is that the operator $\cQ$ needed above can be
explicitly reconstructed (up to a factor on the right) using the functions
$\{\psi_j\}$ from its kernel \eqref{6.5a}. More precisely, if $[m_1,m_2]$ is
the support of $\cQ$, with $m_2-m_1=m$ we can write the operator $\cQ$
in the form
\begin{equation}                        \label{6.9}
\cQ = \fa(\gamma)
\left\lvert\begin{array}{cccl}
\psi_1(\gamma+m_2)&\cdots&\psi_m(\gamma+m_2)& E^{m_2}   \\
\vdots      &      &\vdots      &\vdots     \\
\psi_1(\gamma+k)&\cdots&\psi_m(\gamma+k)& E^{k}     \\
\vdots      &      &\vdots      &\vdots     \\
\psi_1(\gamma+m_1)&\cdots&\psi_m(\gamma+m_1)& E^{m_1}
\end{array}\right\rvert,
\end{equation}
where $\fa(\gamma)$ is some function in $\gamma$, and the determinant is
expanded from left to right, i.e. the shift operators $E^k$ are pulled to
the right. Indeed, it is obvious that the right-hand side of \eqref{6.9}
defines a difference operator with support $[m_1,m_2]$, whose kernel coincides
with the kernel of $\cQ$. But two operators having the same supports and
kernels must differ by a factor, which proves \eqref{6.9}.

Another useful way to describe the Darboux process \eqref{6.1} is to
define a new operator $\cL$ by the following formula
\begin{equation}                        \label{6.10}
\cL=\fh(L_0)=\prod_{i=1}^m(L_0-x_{j_i}\Id).
\end{equation}
Using \eqref{6.1}-\eqref{6.4} one can easily
show by induction that
\begin{equation}                        \label{6.11}
\cL=\cP\cQ,
\end{equation}
where
\begin{equation}                        \label{6.12}
\cP=P_0P_1\cdots P_{m-1}.
\end{equation}
Thus, from \eqref{6.2} we get
\begin{equation}                        \label{6.13}
\hat{\cL}:=\cQ\cP=\prod_{i=1}^m(\hL-x_{j_i}\Id)=\fh(\hL),
\end{equation}
i.e. we have the Darboux map
\begin{equation}                        \label{6.14}
\cL=\fh(L_0)=\cP\cQ\rightarrow \hat{\cL}=\cQ\cP=\fh(\hL).
\end{equation}

\subsection{The main result}
Now we are ready to formulate our main result. We show in \thref{th6.2}
below, that if
\begin{equation}                        \label{6.15}
d=\pm q^{l/2}\text{ and } a=dq^{\alpha},
\text{ for some } l,\alpha\in\N,\; l,\alpha\geq 1,
\end{equation}
then, the second-order difference operator $\hL$ obtained by iteration of
the Darboux process \eqref{6.1}, starting with the Askey-Wilson operator
$L_0=\La{a}$, and iterating at any subset $\{x_{j_1},\dots,x_{j_m}\}$ of the
points
\begin{equation}                        \label{6.16}
x_k=z_k+z_k^{-1},\quad z_k=dq^k=\pm q^{l/2+k}\quad 0\leq k\leq \alpha-1,
\end{equation}
possesses a two-dimensional space of eigenfunctions which
are also eigenfunctions of a $q$-difference operator in the spectral
variable $z$.

Notice that in this case the functions $\{\psi_j\}$ must be linear
combinations of the functions $r_{\gamma}(z_{j_k})$ and $s_{\gamma}(z_{j_k})$,
computed in \prref{pr5.1}, i.e.
\begin{equation}                        \label{6.17}
\psi_k(\gamma)= \mu_k \raz{a}{z_{j_k}}+
        \nu_k \saz{a}{z_{j_k}}, \quad k=1,2,\dots,m,
\end{equation}
for some constants $\{\mu_k,\nu_k\}_{k=1}^{m}$.
Moreover the operator $L_0$ is diagonalizable (in the vector space $\Ker\cQ$),
i.e. $\sigma_k=0$ for all $k$ in \eqref{6.5b}. It is interesting that this
is a purely $q$-phenomenon. In the $q=1$ case, one needs to use generalized
eigenvectors, see \cite{GY} and \cite{HI3}.

\begin{theorem}\label{th6.2}                    

Assume that \eqref{6.15} holds, i.e.
$a=\pm q^{\frac{l}{2}+\alpha}$ and $d=\pm q^{\frac{l}{2}}$,
$l,\alpha\in\N$, $l,\alpha\geq 1$, with the same choice of sign for
$a$ and $d$. The second-order difference operator $\hat{L}$ which is obtained
by iteration of the Darboux transformation \eqref{6.1},
starting with the associated Askey-Wilson operator $L_0=\La{a}$ \eqref{2.1},
and adding bound states at the points
\begin{align}
&\Big\{x_k=dq^k+(dq^k)^{-1}=
\pm\left(q^{\frac{l}{2}+k}+q^{-\left(\frac{l}{2}+k\right)}\right),\nonumber\\
&\qquad\qquad\qquad\qquad
k\in\{j_1,j_2,\dots,j_m\}\subset\{0,1,\dots,\alpha-1\}\Big\},    \label{6.18}
\end{align}
possesses a two-dimensional space of eigenfunctions which are
also eigenfunctions of a $q$-difference operator in the spectral
variable $z$.
\end{theorem}

If $\Psi(\gamma,z)$ satisfies \eqref{2.9} (i.e. it is a linear combination of
$R_{\gamma}$ and $S_{\gamma}$ defined by \eqref{2.8}), then the function
$\hPsi$ defined by
\begin{equation}                        \label{6.19}
\hPsi = \cQ \Psi(\gamma,z),
\end{equation}
satisfies
\begin{equation}                        \label{6.20}
\hL\hPsi =(z+z^{-1})\hPsi.
\end{equation}
Equation \eqref{6.20} follows immediately from \eqref{6.2}. So, it remains
to show that there exists a $q$-difference operator $\hB$ in $z$, such
that
\begin{equation}                        \label{6.21}
\hB\hPsi = \hLambda \hPsi,
\end{equation}
for some function $\hLambda$.

Before we start with the proof of this fact, we shall make two simple but
important remarks.

\begin{remark}\label{re6.3}                     
So far, we have not said anything about the normalization and the supports
of the operators $\{P_k,Q_k\}_{k=0}^{m-1}$, or equivalently, about the
normalization $\fa(\gamma)$ and the support $[m_1,m_2]$ of $\cQ$ in
\eqref{6.9}. The point is that this is not important for the bispectral
relation \eqref{6.21}. Indeed, if \eqref{6.21} holds, and if
\begin{equation}                        \label{6.22}
\cQ'= {\mathfrak{a}}'(\gamma)E^{m'}\cQ\text{ and }
{\hat{\Psi}}'(\gamma,z)=\cQ'\Psi(\gamma,z)=
    {\mathfrak{a}}'(\gamma){\hat{\Psi}}(\gamma+m',z),
\end{equation}
for some function ${\mathfrak{a}}'(\gamma)$ and some $m'\in\Z$, then it
follows from \eqref{6.21} and \eqref{6.22} that
\begin{equation}                        \label{6.23}
\hB{\hat{\Psi}}'(\gamma,z) = \hLambda'{\hat{\Psi}}'(\gamma,z),
\end{equation}
where $\hLambda'={\hat{\Lambda}}_{\gamma+m'}$, i.e. we still have a
$q$-difference equation in $z$. Thus, the choice of the support $[m_1,m_2]$
and of the normalizing function $\fa(\gamma)$ in \eqref{6.9} is not
essential for the bispectrality.
\end{remark}

\begin{remark}\label{re6.4}                     
From \prref{pr4.1} we know that $\La{a}$ is a Darboux transform from
$\La{aq}$ at $a+a^{-1}$. Thus, we can think of $\hL$ also as a Darboux
transformation from the operator $\La{aq}$ at the points
$\{x_{j_1},x_{j_2},\dots,x_{j_m},a+a^{-1}\}$. From \eqref{4.4}, \eqref{4.5},
and \eqref{6.2} one can deduce that
\begin{equation}                        \label{6.24}
\hL\cQ'=\cQ'\La{aq}\text{ with }\cQ'=\cQ\Qa.
\end{equation}
Let us define functions $\psi'_k(\gamma)$ for $k=1,2,\dots,m+1$ by
\begin{equation}                        \label{6.25}
\psi'_k(\gamma)=(1-a/z_{j_k})\mu_k\raz{aq}{z_{j_k}}+\nu_k\saz{aq}{z_{j_k}}
\text{ for }1\leq k\leq m,
\end{equation}
where $\mu_k$ and $\nu_k$ are the constants in \eqref{6.17}, and
\begin{equation}                        \label{6.26}
\psi'_{m+1}(\gamma)=\saz{aq}{a}.
\end{equation}
Using now equations \eqref{4.6a}-\eqref{4.6b}, we see that
\begin{equation}                        \label{6.27}
\Qa\psi'_k(\gamma)=\frac{(z_{j_k}-a)(1-az_{j_k})}{az_{j_k}}\psi_k(\gamma)
\qquad\text{for }1\leq k\leq m,
\end{equation}
which combined with \eqref{6.24} gives
\begin{equation}                        \label{6.28}
\ker\cQ'=\Span\left\{\psi'_1(\gamma),\psi'_2(\gamma),\dots,
    \psi'_{m+1}(\gamma)\right\}.
\end{equation}
Thus, we can think of $\hL$ as a Darboux transformation from $\La{aq}$, and
the kernel of the corresponding intertwining operator is given explicitly
by \eqref{6.25} and \eqref{6.26}.
\end{remark}

\begin{proof}[Proof of \thref{th6.2}]
From \reref{re6.4}, it follows that we can assume $m$ is even, i.e.
$m=2\bm$ for some $\bm\in\N$. \reref{re6.3} tells us that without any
restriction, we can fix the support of $\cQ$ to be $[-\bm,+\bm]$ (i.e.
$m_2=-m_1=\bm$) in \eqref{6.9}, and we can choose an appropriate
normalization factor $\fa(\gamma)$ in \eqref{6.9}. Let us define
$\fa(\gamma)$ as
\begin{equation}                    \label{6.29}
\fa(\gamma):=\left\{
\begin{array}{ll}
1 & \text{if }\bm \text{ is even;}\\
\frac{1}{q^{-\gamma}-abcdq^{\gamma-1}} & \text{if } \bm \text{ is odd.}
\end{array}\right.
\end{equation}
Notice that
\begin{equation}                        \label{6.30}
I(\fa(\gamma))=(-1)^{\bm}\fa(\gamma),
\end{equation}
i.e. up to a sign, $\fa(\gamma)$ is an $I$-invariant function.
From \leref{le5.2}, \eqref{6.17}, and \eqref{6.9} it follows that
$\cQ$ is a difference operator with coefficients, which are rational
functions in $q^{\gamma}$, i.e. $\cQ\in\RE$. From \leref{le5.2} and
\eqref{6.17} we know that the kernel of $\cQ$ consists of $I$-invariant
functions, i.e.
\begin{equation*}
I(\psi_k(\gamma))=\psi_k(\gamma) \text{ for }1\leq k\leq m.
\end{equation*}
From this relation and the definition of the involution $I$, one can
easily check that
\begin{equation}                        \label{6.31}
I(\psi_k(\gamma+i))=\psi_k(\gamma-i) \text{ for any }i\in\Z.
\end{equation}
Formulas \eqref{6.9}, \eqref{6.30}, \eqref{6.31} show that $\cQ$ is an
$I$-invariant operator, i.e. $\cQ\in\dI$. Since $L_0=\La{a}\in\dI$, we see
immediately from \eqref{6.10} that $\cL\in\dI$, which combined with
\eqref{6.11} gives that $\cP\in\dI$. The proof now follows from
\thref{th3.3} and \thref{th3.1}.
\end{proof}

\begin{remark}\label{re6.5}                     
Let us denote
\begin{equation}                        \label{6.32}
h_a(\gamma)=a^{\gamma}(abq^{\gamma},acq^{\gamma},adq^{\gamma};q)_{\infty},
\end{equation}
and define a new second-order difference operator $\fL$ by conjugating
the Askey-Wilson operator $\La{a}$ in \eqref{2.1} with $h_a(\gamma)$, i.e.
\begin{equation}                        \label{6.33}
\fL=h_a(\gamma)^{-1}\,\La{a}\,h_a(\gamma).
\end{equation}
If we write $\fL$ in the form
\begin{equation}                        \label{6.34}
\fL=\fA E + \fB\Id +\fC E^{-1},
\end{equation}
then a straightforward computation shows that the coefficients of $\fL$
are given by the following formulas
\begin{subequations}                    \label{6.35}
\begin{align}
\fA&=\frac{h_a(\gamma+1)}{h_a(\gamma)}A_{\gamma}=
\frac{1-abcdq^{\gamma-1}}{(1-abcdq^{2\gamma-1})(1-abcdq^{2\gamma})}
                                \label{6.35a}\\
\fC&=\frac{h_a(\gamma-1)}{h_a(\gamma)}C_{\gamma}
=(1-q^{\gamma})(1-abq^{\gamma-1})(1-acq^{\gamma-1})(1-adq^{\gamma-1})
                                    \nonumber\\
&\qquad\times \frac{(1-bcq^{\gamma-1})(1-bdq^{\gamma-1})(1-cdq^{\gamma-1})}
{(1-abcdq^{2\gamma-2})(1-abcdq^{2\gamma-1})}            \label{6.35b}\\
\fB&=B_{\gamma}=q^{\gamma-1}\frac{(1+abcdq^{2\gamma-1})(sq+s'abcd)
-q^{\gamma-1}(1+q)abcd(s+s'q)}
{(1-abcdq^{2\gamma-2})(1-abcdq^{2\gamma})},         \label{6.35c}
\end{align}
\end{subequations}
where $s=a+b+c+d$ and $s'=a^{-1}+b^{-1}+c^{-1}+d^{-1}$.
Notice that the operator $\fL$ is symmetric in $a,b,c,d$ (in fact, this was
the operator originally introduced by Askey and Wilson, see \cite[p.~5]{AW}).
Therefore, we can produce bispectral operators as in \thref{th6.2}, if any
two of the parameters $(a,b,c,d)$ are related in the same manner as $a$
and $d$ in \eqref{6.15}.
\end{remark}

\thref{th6.2} and \reref{re6.5} tell us that if two of the parameters
$a,b,c,d$ are related by a formula similar to \eqref{6.15}, we can
construct bispectral operators by adding bound states to the left
or to the right of the continuous spectrum, depending on the sign of these
two parameters (they must have the same signs!). Assume now that we have
two positive parameters, satisfying \eqref{6.15} with $+$,
and the other two parameters are negative and satisfy \eqref{6.15} with
$-$. Then it is possible to construct bispectral operators by adding
bound states both to the left and to the right of the continuous spectrum.
The corollary below explains what modifications are necessary in the
proof of \thref{th6.2} in this case. The operators considered below contain
the so called continuous $q$-Jacobi operators which can be obtained from
the Askey-Wilson's operator by taking
\begin{equation*}
a=q^{\alpha+\frac{1}{2}}, \quad b=-q^{\beta+\frac{1}{2}}, \quad
c=-q^{\frac{1}{2}}, \quad d=q^{\frac{1}{2}}.
\end{equation*}
This substitution was used by Rahman in \cite{R1}. For a different
parametrization of the continuous $q$-Jacobi polynomials, as well as
relations between the different parametrizations see \cite[p.~83]{KS}.

\begin{corollary}\label{co7.6}              
Assume that
\begin{equation}                        \label{6.36}
a=q^{l_1/2+\alpha}, \quad b=-q^{l_2/2+\beta}, \quad
c=-q^{l_2/2}, \quad d=q^{l_1/2},
\end{equation}
where $\alpha,\beta,l_1,l_2$ are positive integers. Let us define
\begin{align*}
J&=\Big\{q^{l_1/2+k}+q^{-l_1/2-k},\;0\leq k\leq \alpha-1 \Big\},\\
K&=\Big\{-\big(q^{l_2/2+k}+q^{-l_2/2-k}\big),\;0\leq k\leq \beta-1\Big\},
\end{align*}
and let $X=\{x_{j_1},x_{j_2},\dots,x_{j_m}\}$ be a subset of $J\cup K$.
The second-order difference operator $\hat{L}$ which is obtained
by iteration of the Darboux transformation \eqref{6.1},
starting with the associated Askey-Wilson operator $L_0=\La{a}$ \eqref{2.1},
and adding bound states at the points of $X$ possesses a two-dimensional
space of eigenfunctions which are also eigenfunctions of a
$q$-difference operator in the spectral variable $z$.
\end{corollary}

\begin{proof} As in the proof of \thref{th6.2}, we can assume that $m$ is
even, i.e. $m=2\bar{m}$ and the support of $\cQ$ is $[-\bar{m},\bar{m}]$.
From \leref{le5.2} we know that the eigenfunctions of $L$, corresponding
to eigenvalues $x_{j_k}\in J$ are $I$-invariant. The main point is to
see what happens with the eigenfunctions $\psi_k(\gamma)$, corresponding to
$x_{j_k}\in K$. Let us denote by $\fL^b$ the Askey-Wilson operator with
parameters $a,b,c,d$ replaced by $b,c,d,a$, respectively, i.e.
\begin{equation}                        \label{6.37}
\fL^b=L_{b,c,d,a}.
\end{equation}
Then using the notations in \reref{re6.5} we will have
\begin{equation}                        \label{6.38}
\fL=h_b(\gamma)^{-1}\,\fL^b\,h_b(\gamma),
\end{equation}
where $h_b(\gamma)$ is analogous to $h_a(\gamma)$ defined in \eqref{6.32}:
\begin{equation}                        \label{6.39}
h_b(\gamma)=b^{\gamma}(baq^{\gamma},bcq^{\gamma},bdq^{\gamma};q)_{\infty}.
\end{equation}
From \eqref{6.33} and \eqref{6.38} it follows that
\begin{equation}                        \label{6.40}
\La{a}=\frac{h_a(\gamma)}{h_b(\gamma)}\,\fL^b\,
\frac{h_b(\gamma)}{h_a(\gamma)}.
\end{equation}
This equation simply means that the correspondence between the eigenfunctions
of $\fL^b$ and $\La{a}$ is given by multiplication by the function
$h_a(\gamma)/h_b(\gamma)$, i.e. if
\begin{equation*}
\Psi_a(\gamma)=\frac{h_a(\gamma)}{h_b(\gamma)}{\Psi}_b(\gamma),
\end{equation*}
then
\begin{equation*}
\La{a}\Psi_a(\gamma)=\lambda\Psi_a(\gamma)\qquad\text{if and only if}\qquad
\fL^b {\Psi}_b(\gamma)=\lambda{\Psi}_b(\gamma).
\end{equation*}
If we take an eigenfunction of $\fL^b$ with eigenvalue $x_{j_k}\in K$,
then this function must be $I$-invariant by \leref{le5.2}, interchanging
the roles of the parameters $(a,d)$ and $(b,c)$. Thus, it remains to see
how $h_a(\gamma)/h_b(\gamma)$ is transformed under the involution $I$.
Using \eqref{6.36}, \eqref{6.32} and \eqref{6.39} we get
\begin{equation}                        \label{6.41}
\frac{h_a(\gamma)}{h_b(\gamma)}
=\Big(\frac{a}{b}\Big)^{\gamma}
\frac{(acq^{\gamma},adq^{\gamma};q)_{\infty}}
{(bcq^{\gamma},bdq^{\gamma};q)_{\infty}}
=(-1)^{\gamma}h(\gamma),
\end{equation}
where
\begin{equation}                        \label{6.42}
h(\gamma)=q^{[(l_1-l_2)/2+\alpha-\beta]\gamma}\,
\frac{(q^{\gamma+\alpha+l_1},-q^{\gamma+(l_1+l_2)/2+\alpha};q)_{\infty}}
{(q^{\gamma+\beta+l_2},-q^{\gamma+(l_1+l_2)/2+\beta};q)_{\infty}}.
\end{equation}
Using the definitions at the beginning of Section 2 we can write $h(\gamma)$
in the form
\begin{equation}                        \label{6.43}
h(\gamma)=
\frac{q^{[(l_1-l_2)/2+\alpha-\beta]\gamma}}
{(q^{\gamma+\beta+l_2};q)_{\alpha-\beta+l_1-l_2}\,
(-q^{\gamma+(l_1+l_2)/2+\beta};q)_{\alpha-\beta}},
\end{equation}
where
\begin{equation*}
(y;q)_n=\frac{1}{(yq^n;q)_{-n}}\text{ for }n<0.
\end{equation*}
Notice that if $l_1-l_2$ is even, then $h(\gamma)$ is a rational function
in $q^{\gamma}$, i.e. $h(\gamma)\in\RE$. However, for $l_1-l_2$ odd,
$h(\gamma)$ has an extra factor $q^{\gamma/2}$, i.e.
$h(\gamma)=q^{\gamma/2}\times$(rational function in $q^{\gamma}$).

In the case considered here (i.e. assuming \eqref{6.36}), the involution
$I$ acts on rational functions in $q^{\gamma}$ by the following formula
\begin{equation*}
I(q^{\gamma})=q^{-\gamma-\alpha-\beta-l_1-l_2+1}.
\end{equation*}
We can naturally extend the involution $I$ to an involution acting on
$\mathcal{R}\{q^{\gamma/2},E\}$ by defining
\begin{equation}                        \label{6.44}
I(q^{\gamma/2})=q^{-(\gamma+\alpha+\beta+l_1+l_2-1)/2}.
\end{equation}
Using \eqref{6.43} and \eqref{6.44} one can check that
\begin{equation}                        \label{6.45}
I(h(\gamma))=(-1)^{\alpha+\beta+l_1+l_2}h(\gamma),
\end{equation}
i.e. up to a sign the function $h(\gamma)$ is $I$-invariant.
Finally, let us denote by $m'$ and $m''$ the number of
elements in the sets $X\cap J$ and $X\cap K$, respectively. The
computations above can be summarized as follows:
\begin{itemize}
\item $m'$ of the functions $\psi_k$, generating the kernel of $\cQ$,
are eigenfunctions of $\La{a}$ corresponding to eigenvalues $x_{j_k}\in J$ and
these functions are $I$-invariant rational functions in $q^{\gamma}$;
\item $m''$ of the functions $\psi_k$, generating the kernel of $\cQ$,
are eigenfunctions of $\La{a}$ corresponding to eigenvalues $x_{j_k}\in K$ and
these functions are of the form
\begin{equation*}
\psi_k(\gamma)=(-1)^{\gamma}\xi_k(\gamma),
\end{equation*}
where
\begin{equation}                        \label{6.46}
I(\xi_k(\gamma))=(-1)^{\alpha+\beta+l_1+l_2}\xi_k(\gamma),
\end{equation}
and
\begin{equation}                        \label{6.47}
\xi_k(\gamma)=
\left\{ \begin{array}{ll}
\text{rational function in }q^{\gamma} & \text{ if }l_1-l_2 \text{ is even;}\\
q^{\gamma/2}\times(\text{rational function in }q^{\gamma})
& \text{ if }l_1-l_2 \text{ is odd.}
\end{array}\right.
\end{equation}
\end{itemize}
Without any restriction, we can assume that $x_{j_k}\in J$ for
$k=1,2,\dots,m'$ and $x_{j_s}\in K$ for
$s=m'+1,\dots,m''$. Let us denote by $(a_{ij};b_{ij};c_i)$ the
$(m+1)\times(m+1)$ matrix with rows numbered from $\bar{m}$ to $-\bar{m}$
(i.e. $i=\bar{m},\bar{m}-1,\dots,-\bar{m}$) and columns from 1 to $m+1$,
with entries $a_{ij}$ for the first $m'$ columns, entries $b_{ij}$
for the next $m''$ columns, and entries $c_i$ in the last column.
Then, from \eqref{6.9} we see that the operator $\cQ$ has the form
\begin{equation}                                    \label{6.48}
\cQ = \fa(\gamma) (-1)^{m''\gamma}
\det(\psi_j(\gamma+i);(-1)^i\xi_j(\gamma+i);E^{i}),
\end{equation}
where $\fa(\gamma)$ is an appropriate function. From \eqref{6.46}
it follows that
\begin{equation*}
I(\xi_j(\gamma+i))=(-1)^{\alpha+\beta+l_1+l_2}\xi_j(\gamma-i),
\end{equation*}
which combined with \eqref{6.48} shows that
\begin{equation}                        \label{6.49}
I(\cQ)=(-1)^M\, I(\fa(\gamma) (-1)^{m''\gamma})
\det(\psi_j(\gamma+i);(-1)^i\xi_j(\gamma+i);E^{i}),
\end{equation}
where $M=\bar{m}+(\alpha+\beta+l_1+l_2)m''$. From \eqref{6.47} and
\eqref{6.48} we see that
\begin{equation}                        \label{6.50}
\cQ=\fa(\gamma) (-1)^{m''\gamma}\times
\left\{\begin{array}{ll}
\text{operator in }\RE & \text{if }m''(l_1-l_2) \text{ is even;}\\
q^{\gamma/2}\times(\text{operator in }\RE) &
\text{if }m''(l_1-l_2) \text{ is odd.}
\end{array}\right.
\end{equation}
Let us now define $\fa(\gamma)$ as follows
\begin{equation}                        \label{6.51}
\fa(\gamma)=\frac{1}{(-1)^{m''\gamma}(q^{\epsilon_1\gamma/2}+
\epsilon_2q^{-\epsilon_1(\gamma+\alpha+\beta+l_1+l_2-1)/2})},
\end{equation}
where
\begin{equation*}
\epsilon_1=\left\{\begin{array}{ll}
2 & \text{ if }m''(l_1-l_2) \text{ is even;}\\
1 & \text{ if }m''(l_1-l_2) \text{ is odd;}
\end{array}\right.
\end{equation*}
and
\begin{equation*}
\epsilon_2=(-1)^M=
\left\{\begin{array}{rl}
 1 & \text{ if } M \text{ is even;}\\
-1 & \text{ if } M \text{ is odd.}
\end{array}\right.
\end{equation*}
From \eqref{6.48}, \eqref{6.49}, \eqref{6.50} it follows that
with this definition of $\fa(\gamma)$, $\cQ$ becomes an $I$-invariant
operator from $\RE$, and the proof continues as in \thref{th6.2}.
\end{proof}

\section{A $q$-analogue of the Krall-Jacobi polynomials, with a bound
state off the continuous spectrum}\label{se7}

In this section we illustrate \thref{th6.2} on the special case
\begin{equation}\label{7.1}
a=q^{\frac{3}{2}},\quad d=q^{\frac{1}{2}},\quad b\;
\mbox{and}\;c\quad\mbox{arbitrary},
\end{equation}
adding a bound state off the continuous spectrum of $\La{a}$ at the
point
\begin{equation}\label{7.2}
x_{0}=q^{\frac{1}{2}}+q^{-\frac{1}{2}}.
\end{equation}
This corresponds to the choice
\begin{equation}\label{7.3}
l=1,\quad \alpha=1,\quad k=0.
\end{equation}
\subsection{Illustrating the general theory}
According to \reref{re6.4}, this case must be handled by first expressing
$\La{a}$ as a Darboux transform (without free parameter) from $\La{aq}$ at
$a+a^{-1}$, and then performing a Darboux
transform (with free parameter) of $\La{a}$ at $x_{0}$. Thus we
perform the chain of elementary Darboux transformations
\begin{multline}\label{7.4}
\La{aq}=(a+a^{-1})\;\Id+\Pa\Qa\curvearrowright\\
\La{a}=(a+a^{-1})\;\Id+\Qa\Pa
=x_{0}\;\Id+PQ\\
\curvearrowright \hL=x_{0}\;\Id+QP.
\end{multline}

We now explain how Theorems 3.1 and 3.3 (which form the core of the
general proof of \thref{th6.2}) apply to this
example. We can factorize
\begin{equation}\label{7.5}
\cL=\big(\La{aq}-x_{0}\Id\big)\big(\La{aq}-(a+a^{-1})\Id\big),
\end{equation}
as
\begin{equation}\label{7.6}
\cL=\cP\cQ\quad \mbox{with}\quad
\cP=\Pa P\quad \mbox{and}\quad \cQ=Q\Qa.
\end{equation}
The explicit formulas for $\Pa$ and $\Qa$ are given in \eqref{4.2}
and \eqref{4.3} respectively. The kernel $\psi_{\gamma}$ of the operator
$\La{a}-x_{0}\;\Id$ is readily obtained by substitution of \eqref{7.1} and
\eqref{7.3} in \eqref{5.2} and \eqref{5.3} (using \eqref{5.1}, \eqref{5.4},
\eqref{5.6} and \eqref{5.7}), which leads to
\begin{equation}\label{7.7}
\psi_{\gamma}=\frac{q^{\gamma}}{(1-bq^{\gamma+\frac{1}{2}})
(1-cq^{\gamma+\frac{1}{2}})}\Bigg(\mu+\frac{q^{\gamma}}{(1-q^{\gamma+1})
(1-bcq^{\gamma})}\Bigg),
\end{equation}
with $\mu$ a free parameter. This gives
\begin{equation}\label{7.8}
Q=\varphi_{\gamma}(\psi_{\gamma}E-\psi_{\gamma+1}),\quad
P=-\frac{\Ca{a}}{\psi_{\gamma}\;\varphi_{\gamma-1}}E^{-1}+
\frac{\Aa{a}}{\psi_{\gamma}\;\varphi_{\gamma}}.
\end{equation}
In this last formula $\varphi_{\gamma}$ is as in \eqref{4.1} and corresponds
to the factor $\fa(\gamma)$ defined in \eqref{6.29} (with
$\bm=1$ and $a$ replaced by $aq$), and $\Aa{a}$ and $\Ca{a}$ are defined as
in \eqref{2.2a} and \eqref{2.2b}; it is understood that the values
\eqref{7.1} for $a$ and $d$ are to be substituted everywhere.

The following identities are equivalent to \eqref{4.5}:
\begin{equation*}
\frac{\varphi_{\gamma}}{\varphi_{\gamma+1}}\Ag{\gamma+1}=\Aa{aq}\quad
\mbox{and} \quad
\frac{\varphi_{\gamma}}{\varphi_{\gamma-1}}\Cg{\gamma}=\Ca{aq}.
\end{equation*}
Using these identities, by a straightforward
computation, we obtain that
\begin{multline}\label{7.9}
\cQ=Q\Qa=\Aa{aq}\;\psi_{\gamma}\;E\\
-(\Aa{a}\;\psi_{\gamma+1}+\Cg{\gamma+1}\;\psi_{\gamma})\;\Id+\Ca{aq}\;
\psi_{\gamma+1}\;E^{-1},
\end{multline}
and
\begin{multline}\label{7.10}
\cP=\Pa P=\frac{\Aa{aq}}{\psi_{\gamma+1}}E\\
-\Bigg(\frac{\Ag{\gamma}}{\psi_{\gamma}}+\frac{\Cg{\gamma+1}}
{\psi_{\gamma+1}}\Bigg)\;\Id
+\frac{\Ca{aq}}{\psi_{\gamma}}\;E^{-1}.
\end{multline}

Since $\cL$ in \eqref{7.5} is a constant coefficient polynomial in
$\La{aq}$, the relevant involution $I$
is obtained by replacing $a$ with $aq$ in \eqref{3.15}, that is
\begin{equation*}
I(q^{\gamma})=\frac{1}{abcdq^{\gamma}}=\frac{1}{bcq^{\gamma+2}}\quad
\mbox{and}\quad I(E)=E^{-1}.
\end{equation*}
We know from \eqref{3.18} that $I(\Aa{aq})=\Ca{aq}$. Using \eqref{2.2a},
\eqref{2.2b} and \eqref{7.7}, it is straightforward to check that
\begin{equation}\label{7.11}
I(\Ag{\gamma})=\Cg{\gamma+1}\quad \mbox{and}\quad
I(\psi_{\gamma})=\psi_{\gamma+1}.
\end{equation}
This makes clear that the operators $\cQ$ and $\cP$ as
expressed in \eqref{7.9} and \eqref{7.10} are $I$-invariant.
To alleviate the notations, we shall write $L$ and $\Lambda$
to denote the following difference operators:
\begin{equation*}
L\equiv \La{aq},\quad \Lambda\equiv
q^{-\gamma}(1-q^{\gamma})(1-abcdq^{\gamma})\;\Id,\quad
\mbox{with}\;a,b,c,d\;\mbox{as in}\;\eqref{7.1}.
\end{equation*}
According to \thref{th3.3}, we must have that
\begin{equation*}
\cP=U\Gamma^{-1}\quad \mbox{and}\quad \cQ=\Theta^{-1}V,
\end{equation*}
with $U,V\in\cB=\langle L,\Lambda\rangle$ and $\Gamma, \Theta\in
\cK=\langle\Lambda\rangle$.

Using the symbolic software MACSYMA, we have obtained the
following explicit formulas for $\Gamma, \Theta, U$ and $V$. It
follows from \eqref{7.11} that the function
$\psi_{\gamma}\psi_{\gamma+1}$ is $I$-invariant, hence it can be
written as a rational function of $\Lambda$. Explicitly:
\begin{gather*}
\psi_{\gamma}\psi_{\gamma+1}\;\Id=\Gamma(\Lambda)\Theta^{-1}(\Lambda),
\quad \mbox{with}\\
\Gamma(\Lambda)=\big[(q-1)(bc-1)\mu+1\big]\big[(q^2-1)(bcq-1)\mu+q\big]\Id\\
+\mu\big[q+1+\mu(q-1)(2bcq-q+bc-2)\big]\Lambda+\mu^2\Lambda^2,\\
\Theta(\Lambda)=q^{-1}\big[\Lambda+(q^2-1)(bc-1)\Id\big]
\big[\Lambda+(q-1)(bcq-1)\Id\big]\\
\times\big[\Lambda+(bq^{1/2}-1)(cq^{3/2}-1)\Id\big]
\big[\Lambda+(bq^{3/2}-1)(cq^{1/2}-1)\Id\big].
\end{gather*}
Then, one finds that
\begin{align*}
U=\;&(r_{1}\;\Lambda+r_{2}\;\Lambda^{2}+r_{3}\;\Lambda^{3})\;L+L\;(s_{1}\;
\Lambda+s_{2}\;
\Lambda^{2}+s_{3}\;\Lambda^{3})\\&+t\;L+u_{0}\;\Id+u_{1}\;
\Lambda+u_{2}\;\Lambda^{2}+u_{3}\;\Lambda^{3},\\
V=\;&(v_{1}\;\Lambda+v_{2}\;\Lambda^{2}+v_{3}\;\Lambda^{3})\;
L+L\;(w_{1}\;\Lambda+w_{2}\;
\Lambda^{2}+w_{3}\;\Lambda^{3})\\&+x\;\Lambda\; L\;\Lambda+y_{0}\;\Id+
y_{1}\;\Lambda+y_{2}\;\Lambda^{2}+y_{3}\;\Lambda^{3},
\end{align*}
for some appropriate choice of the constants
$r_{i},s_{i},t,u_{i},v_{i},w_{i},x,y_{i}$. The explicit
expressions for these constants are rather complicated and there is
no point in displaying them. \thref{th3.1} (see especially
\eqref{3.10} and \eqref{3.11}) provides the following
explicit bispectral operator
\begin{equation}\label{7.12}
\hat{B}=\fb(V)\fb(U)f^{-1},
\end{equation}
with
\begin{align*}
f=\;&\fb(\cL)=\big(z+z^{-1}-(q^{\frac{1}{2}}+q^{-\frac{1}{2}})\big)
\big(z+z^{-1}-(q^{\frac{3}{2}}+q^{-\frac{3}{2}})\big),\\
\fb(U)=\;&(z+z^{-1})(r_{1}B+r_{2}B^{2}+r_{3}B^{3})+
(s_{1}B+s_{2}B^{2}+s_{3}B^{3})(z+z^{-1})\\
&+t(z+z^{-1})+u_{0}+u_{1}B+u_{2}B^{2}+u_{3}B^{3},\\
\fb(V)=\;&(z+z^{-1})(v_{1}B+v_{2}B^{2}+v_{3}B^{3})+
(w_{1}B+w_{2}B^{2}+w_{3}B^{3})(z+z^{-1})\\
&+xB(z+z^{-1})B+y_{0}+y_{1}B+y_{2}B^{2}+y_{3}B^{3},
\end{align*}
where $B=B_{q^{5/2},b,c,q^{1/2}}(z,D_{z})$ and $\Bz$ is defined as in
\eqref{2.6}. Since $B$ is a second-order $q$-difference operator, $\hat{B}$
in \eqref{7.12} is a $q$-difference operator of order 12.
\subsection{A bispectral operator of order 4}
It is possible for the example discussed above, with $a,b,c,d$ as
in \eqref{7.1}, to produce a bispectral operator of order
4, by performing a bispectral Darboux transformation of the
operator
\begin{equation*}
\mathcal{L}=\big(L_{q^{1/2},b,c,q^{1/2}}-(q^{1/2}+q^{-1/2})\Id\big)^{2},
\end{equation*}
instead of the choice we made in \eqref{7.5}. The idea is to first
obtain the operator $L_{q^{3/2},b,c,q^{1/2}}$ as a (contiguous) Darboux
transform of $L_{q^{1/2},b,c,q^{1/2}}$ at $q^{1/2}+q^{-1/2}$ using
\eqref{4.4} and \eqref{4.5}, and then to perform a (general) upper-lower
Darboux transform (with a free parameter) of $L_{q^{3/2},b,c,q^{1/2}}$
at $q^{1/2}+q^{-1/2}$. In the case $q=1$, such a strategy has already
been explained with all details in
Section 4 of \cite{H}. Since we don't know yet how to produce
a bispectral operator of the lowest possible order for the general situation
described in \thref{th6.2}, we shall only
display the final result. Let us define
\begin{gather}
D_{-2}(z)=\frac{\mu\;q^{3/2}(z-b)(z-c)(z-bq)(z-cq)(z-q^{5/2})(q^{1/2}z-1)}
{(z+q^{1/2})(z+q^{3/2})(z^{2}-1)(z^{2}-q^{2})},\nonumber\\
D_{-1}(z)=q(q+1)(z-b)(z-c)(z-q^{1/2})\nonumber\\
\times \frac{(z+q^{3/2})(q^{1/2}z+1)-\mu\;q^{1/2}
\Big[(q+bc)(z^{2}+q)-(b+c)(q^{2}+1)z\Big]}
{(z+q^{1/2})(z+q^{3/2})(q^{1/2}z+1)(z^{2}-1)},\nonumber\\
D_{1}(z)=D_{-1}(1/z),\qquad D_{2}(z)=D_{-2}(1/z),\nonumber\\
D_{0}(z)=-D_{-2}(z)-D_{-1}(z)-D_{1}(z)-D_{2}(z),\label{7.13}
\end{gather}
and
\begin{equation*}
\hat{\Lambda}(\gamma)=\frac{(q^{\gamma+1}-1)(bcq^{\gamma+1}-1)
\Big[q^{\gamma}(q+1)
+\mu\;(bcq^{\gamma}-1)(q^{\gamma+2}-1)\Big]}{q^{2\gamma}}.
\end{equation*}
Then, the functions
\begin{gather*}
\hat{\Psi}(\gamma,z)=Q\Psi(\gamma,z),\quad \mbox{with}\\
\Psi(\gamma,z)=R_{\gamma}(q^{3/2},b,c,q^{1/2};z)\quad\mbox{or}
\quad S_{\gamma}(q^{3/2},b,c,q^{1/2};z),
\end{gather*}
and with $Q$ as in \eqref{7.8}, $R_{\gamma}(a,b,c,d;z)$ and
$S_{\gamma}(a,b,c,d;z)$ as in \eqref{2.8a} and \eqref{2.8b},
satisfy the pair of equations
\begin{gather*}
\hat{L}\hat{\Psi}(\gamma,z)=(z+1/z)\hat{\Psi}(\gamma,z),\quad \hat{L}\;
\mbox{as in}\;\eqref{7.4},\\
\sum_{i=-2}^{2}D_{i}(z)\hat{\Psi}(\gamma,q^{i}z)=\hat{\Lambda}(\gamma)
\hat{\Psi}(\gamma,z).
\end{gather*}

Let us consider the special case of \eqref{7.1} given by
\begin{equation}\label{7.14}
a=q^{3/2},\quad d=q^{1/2},\quad b=-q^{1/2},\quad
c=-q^{\beta+1/2},\;\beta>-1,
\end{equation}
and let us also specialize $\gamma$ to be a positive integer. Then, the
functions
\begin{equation*}
R_{n}(q^{3/2},-q^{1/2},-q^{\beta+1/2},q^{1/2};z), \quad n=0,1,2,\dots,
\end{equation*}
become (up to a normalizing factor) the continuous $q$-Jacobi polynomials
with parameters $(\alpha=1,\;\beta)$ as defined by Rahman (see
\cite{KS}, formula (3.10.14)). When $q\to 1$, these polynomials reduce to the
standard Jacobi polynomials $P_{n}^{(1,\beta)}(x),\;x=(z+z^{-1})/2$
(see \cite{KS}, formula (5.10.2)), with weight function
$(1-x)(1+x)^{\beta}$ on $[-1,1]$.
The new functions
\begin{gather*}
\hat{P}_{0}(z)=1,\\
\hat{P}_{n}(z)\equiv(QR)_{n-1}=\phi_{n-1}
\psi_{n-1}R_{n}(q^{3/2},-q^{1/2},-q^{\beta+1/2},q^{1/2};z)\\
-\phi_{n-1}\psi_{n}R_{n-1}(q^{3/2},-q^{1/2},-q^{\beta+1/2},q^{1/2};z),
\quad n=1,2,3,\dots,
\end{gather*}
obtained from the Darboux process, with $Q$ as in \eqref{7.8},
are (Laurent) polynomials of degree $n$ in the variable
$z+z^{-1}$. These polynomials are eigenfunctions of the
fourth-order $q$-difference operator with coefficients as in \eqref{7.13},
with $b$ and $c$ as in \eqref{7.14}.

Putting the free parameter $\mu$ of the Darboux transform in \eqref{7.7} to be
\begin{equation*}
\mu=\frac{\rho}{2^{\beta+1}(q-1)^2},
\end{equation*}
one can check that, taking the limit $q\to 1$, the polynomials
\begin{equation*}
K_{n}(x)=\lim_{q\to 1}\;(q-1)^{3}\hat{P}_{n}(z),\quad
x=\frac{z+z^{-1}}{2},
\end{equation*}
are the so-called Krall-Jacobi polynomials, which are orthogonal on the
interval $[-1,1]$ for the weight function
\begin{equation}\label{7.15}
(1+x)^{\beta}+\rho\;\delta(x-1),
\end{equation}
where $\delta(x)$ denotes Dirac's delta function. This is one of the family
of orthogonal polynomials satisfying a fourth-order differential
equation that was discovered by H.L. Krall in \cite{Kr2}. Thus, the example
discussed in this section can be viewed as a $q$-deformation of
the Krall-Jacobi polynomials, with a bound state of arbitrary weight
$\rho$ off the continuous spectrum, precisely at the point $x\equiv
(z+z^{-1})/2=(q^{1/2}+q^{-1/2})/2$. When $q\to 1$, this bound state tends to
the boundary $x=1$ of the continuous spectrum of
the Jacobi polynomials, in agreement with \eqref{7.15}.

\section{Askey-Wilson type solitons}\label{se8}

In this section we construct a large family of second-order
difference operators which provide solutions to the Askey-Wilson bispectral
problem as defined in \eqref{1.1} and \eqref{1.2}, within the context of
the theory of commutative rank 1 rings of difference operators.
We start from the simple observation that when $a=-b=1$ and
$c=-d=\sqrt{q}$, the Askey-Wilson operator \eqref{2.1} reduces to
\begin{equation}                        \label{8.1}
\La{a}=E+E^{-1}.
\end{equation}
More precisely, we shall prove the next theorem:
\begin{theorem}\label{th8.1}                
Let $k_{i}, 1\leq i \leq g$, denote $g$ positive
integers. Let us partition these integers into two arbitrary sets
$J$ and $K$ of $\alpha$ and $\beta$ elements respectively,
$\alpha, \beta \geq 0, \alpha+\beta=g$. We build accordingly the
$\alpha$ points $q^{k_{i}/2}+q^{-k_{i}/2}, i \in J$, and the
$\beta$ other points $-(q^{k_{i}/2}+q^{-k_{i}/2}), i \in K$,
and we assume that all these points are distinct.

i) The tridiagonal operator $\hat{L}$ which is obtained by iteration of
the Darboux process starting at $L=E+E^{-1}$, and adding bound
states at these $g=\alpha+\beta$ distinct points, is part of a
commutative rank 1 ring of difference operators, with
spectral curve as in \eqref{1.6}.

ii) The (reduced) wave function $\Psi(n,z)$ of this rank 1 commutative
ring of difference operators, is also a common eigenfunction of a rank 1
commutative ring of $q$-difference operators in the spectral variable z.
\end{theorem}

The proof of this result is closely connected with the idea of "duality"
that we exploited to establish \thref{th2.1}. In the context of rank
1 commutative rings of differential operators, "duality" was first
systematically exploited by G. Wilson \cite{Wi1}.
Before proving \thref{th8.1}, we need to give a
brief summary of the theory of rank 1 commutative rings of
difference operators. We limit ourselves to the case when the
spectrum of the ring is a rational irreducible curve, which is the
only part of the theory that we shall use. The presentation
follows closely our previous work \cite{HI2}.


\subsection{Commutative rings of difference operators and
rational curves}
We denote by $\mathcal{P}=\mathbb{C}[z]$ the space of polynomials
in the variable $z$, and we denote by $\mathcal{R}=\mathbb{C}(z)$
the space of rational functions of $z$.
\begin{definition}\label{de8.2}             
A flag $\mathcal{V}:\cdots\subset
V_{n+1}\subset V_{n}\subset V_{n-1}\subset \cdots$ of subspaces
$V_{n}\subset \mathcal{R}$ is called a rational flag if and only
if there are polynomials $p(z)$ and $r(z)$ (independent of n)
such that
\begin{equation*}
z^{n}p(z)\mathcal{P}\subset V_{n}\subset
z^{n}r^{-1}(z)\mathcal{P},
\end{equation*}
and the codimension of $V_{n}$ in $z^{n}r^{-1}(z)\mathcal{P}$ is
equal to the degree of $r(z)$.
\end{definition}

There is a one-to-one correspondence between rational flags and
affine irreducible rational curves which complete by adding
\emph{two} non-singular points at infinity. The curve is recovered
from the flag as the spectrum of the ring
\begin{multline} \label{8.2}
A_{\mathcal{V}}=\{f(z)\in\mathcal{R}\;\mbox{with poles only
at}\;z=0\;\mbox{and}\;z=\infty,\\
\mbox{such that}\;\exists k\in\mathbb{Z}\;
\mbox{for which}\;f(z)V_{n}\subset V_{n+k}, \forall n\}.
\end{multline}
When the curve is rationally parametrized by $z$, the two non-singular
distinguished points at infinity $P_{\infty}$ and $Q_{\infty}$
correspond respectively to $z=\infty$ and $z=0$.

By definition, the \emph{rank} of a commutative ring of difference
operators is the greatest common divisor of the orders of all the
operators in the ring. There is a dictionary between the affine
rings of irreducible algebraic curves which complete by adding
\emph{two} non-singular points at infinity and the rank 1
commutative rings of difference operators. We refer the reader to
\cite{Mu} for a full account of the theory, as well as for
adequate references. We limit ourselves
to the case of irreducible rational curves as described above.
The general case proceeds along similar lines. To each $V_{n}$ there
corresponds a tau function $\tau_{V_{n}}(t)
\equiv \tau(n,t)$ in the sense of Segal and Wilson \cite{SW},
depending on infinitely many variables $t=(t_{1},t_{2},t_{3},\dots)$.
The \emph{wave} and \emph{adjoint wave functions}, denoted
respectively by $w(n,t,z)$ and $w^{*}(n,t,z)$, are defined by
\begin{align}
w(n,t,z)
&=z^{n}\mbox{exp}(t,z)\frac{\tau(n,t-[z^{-1}])}{\tau(n,t)},\label{8.3}\\
w^{*}(n,t,z)
&=z^{-n}\mbox{exp}^{-1}(t,z)\frac{\tau(n,t+[z^{-1}])}{\tau(n,t)},
                          \text{ with }\label{8.4}\\
\exp(t,z)&=\exp(t_{1}z+t_{2}z^{2}+t_{3}z^{3}+\cdots)
\text{ and }[z]=(z,z^{2}/2,z^{3}/3,\dots).               \nonumber
\end{align}

The plane $V_{n}$ can be recovered from the wave
function as
\begin{equation*}
V_{n}=\mbox{Span}\{w(n,0,z),w(n+1,0,z),w(n+2,0,z),\dots\}.
\end{equation*}
The \emph{dual flag} of subspaces $\mathcal{V}^{*}:\cdots \supset
V_{n+1}^{*}\supset V_{n}^{*}\supset V_{n-1}^{*}\supset\cdots$,
with
\begin{equation*}
V_{n}^{*}=\Span\{w^{*}(n,0,z),w^{*}(n-1,0,z),w^{*}(n-2,0,z),\dots\},
\end{equation*}
defines the same tau function, up to a change of sign,
$\tau_{V_{n}^{*}}(t)=\tau_{V_{n}}(-t)$.

\begin{theorem}\label{th8.3}                    
\emph{(see \cite{HI2}, \cite{Mu} and references therein)}.
The wave function $w(n,t,z)$ is the common eigenfunction of a commutative
rank 1 ring of difference operators which is isomorphic to
$A_{\mathcal{V}}$, as defined in \eqref{8.2}. More precisely, for any
$f\in A_{\mathcal{V}}$, there is a difference operator $L_{f}$ such that
\begin{equation*}
L_{f}w(n,t,z)=f(z)w(n,t,z).
\end{equation*}
If $f$ has a pole of order $i$ at $P_{\infty}$ and a pole of order
$j$ at $Q_{\infty}$, the operator $L_{f}$, thought of as a finite
band matrix, has $i$ diagonals above the main diagonal and $j$
diagonals below it.
\end{theorem}

We illustrate the concepts above on the example of the
simplest rational singular curves, namely those which
are obtained by identifying $2g$ distinct points on the Riemann
sphere $\mathbb{P}^{1}(\mathbb{C})$ in pairs
$\{\lambda_{i},\mu_{i}\},i=1,\dots,g$. We assume that all these
points are distinct from the two distinguished points $Q_{\infty}=0$ and
$P_{\infty}=\infty$. We pick $g$ nonzero arbitrary complex numbers
$\delta_{i}\in\mathbb{C}^{*},1\leq i\leq g$ (i.e. a divisor
of degree $g$ on the curve), and we define
\begin{multline}                    \label{8.5}
V_{n}=\frac{1}{\prod_{i=1}^{g}(z-\lambda_{i})}\Big\{\text{meromorphic
functions $f$ on $\mathbb{P}^{1}(\mathbb{C})$ such that}\\
(f)-n Q_{\infty}\geq 0\text{ on }\mathbb{P}^{1}(\mathbb{C})\setminus
\{\infty\}\text{ and }f(\lambda_{i})=\delta_{i}f(\mu_{i}), \;
i=1,\dots,g\Big\}.
\end{multline}
By the notation $(f)-nQ_{\infty}\geq 0$, we mean that $f$ has a
zero at least of order $n$ at $Q_{\infty}$ if $n\geq 0$ and a pole
of order at most $-n$ at $Q_{\infty}$ if $n<0$. The next lemma
shows that this definition fits within \deref{de8.2}. In order to
establish it, it is useful to introduce the following definitions.
A function $f\in\mathcal{R}$ which admits an expansion around
$z=\infty$ of the form
\begin{equation}                        \label{8.6}
f=c_{n}z^{n}+c_{n-1}z^{n-1}+\cdots,\quad c_{n}\neq 0,
\end{equation}
will be called an \emph{element of order n}. It will also be
useful to consider the non-degenerate bilinear form on
$\mathcal{R}$
\begin{equation}                        \label{8.7}
B(f,g)=\mbox{res}_{z=\infty}f(z)g(z)\mbox{d}z,\quad f,g\in\mathcal{R}.
\end{equation}
\begin{lemma}\label{le8.4}                      
The plane $V_{n}$ defined in \eqref{8.5} satisfies
\begin{equation*}
z^{n}\prod_{i=1}^{g}(z-\mu_{i})\mathcal{P}\subset V_{n}\subset
z^{n}\prod_{i=1}^{g}(z-\lambda_{i})^{-1}\mathcal{P},
\end{equation*}
and the codimension of $V_{n}$ in $z^{n}\prod_{i=1}^{g}
(z-\lambda_{i})^{-1}\mathcal{P}$ is equal to $g$.
\end{lemma}

\begin{proof} The second inclusion as well as the assertion about
the codimension both follow immediately from the definition of
$V_{n}$ in \eqref{8.5}. In order to establish the first inclusion, we
introduce the dual flag of subspaces
\begin{multline}                        \label{8.8}
V_{n}^{*}=\frac{1}{\prod_{i=1}^{g}(z-\mu_{i})}\Big\{\text{meromorphic
functions $f$ on $\mathbb{P}^{1}(\mathbb{C})$ such that}\\
(f)+nQ_{\infty}\geq 0\text{ on }\mathbb{P}^{1}(\mathbb{C})\setminus\{\infty\}
\text{ and }f(\mu_{i})=\delta_{i}^{*}f(\lambda_{i}),\,i=1,\dots,g\Big\},
\end{multline}
with
\begin{equation}                        \label{8.9}
\delta_{i}^{*}=\delta_{i}
\prod_{\begin{subarray}{c} j=1\\ j\neq i\end{subarray}}^{g}
\frac{(\mu_{i}-\mu_{j})(\mu_{i}-\lambda_{j})}
{(\lambda_{i}-\mu_{j})(\lambda_{i}-\lambda_{j})}.
\end{equation}
It follows immediately from \eqref{8.5}, \eqref{8.8} and \eqref{8.9} that
\begin{equation*}
\{\mbox{res}_{z=\lambda_{i}}+\mbox{res}_{z=\mu_{i}}\}w(z)w^{*}(z)\mbox{d}z=0,
\quad \forall \;w(z)\in V_{n}\;\mbox{and}\;\forall \;w^{*}(z)\in V_{m}^{*}.
\end{equation*}
For $w(z)\in V_{n}$ and $w^{*}(z)\in V_{m}^{*}$ with $n\geq m$,
$w(z)w^{*}(z)$ has no pole at $z=0$. Thus, by the residue
theorem, we deduce that
\begin{multline*}
\mbox{res}_{z=\infty}w(z)w^{*}(z)\mbox{d}z=
-\sum_{i=1}^{g}\{\mbox{res}_{z=\lambda_{i}}+
\mbox{res}_{z=\mu_{i}}\}w(z)w^{*}(z)\mbox{d}z=0,\\
\forall \;w(z)\in V_{n}\;\mbox{and}\;\forall\; w^{*}(z)\in
V_{m}^{*}\;\mbox{with}\;n\geq m,
\end{multline*}
showing that $V_{n}^{*}\subseteq V_{n}^{\perp}$, with
$V_{n}^{\perp}$ the orthogonal of $V_{n}$ with respect to the
bilinear form $B$ defined in \eqref{8.7}. Since $V_{n}$ has a basis
of elements of orders $n,n+1,n+2,\dots$ (see \eqref{8.6} for the
definition), and $V_{n}^{*}$ has a
basis of elements of orders $-n,-n+1,-n+2,\dots$, we deduce that
$V_{n}^{*}=V_{n}^{\perp}$. By the definition \eqref{8.8} of
$V_{n}^{*}$, it implies that
\begin{equation*}
V_{n}^{\perp}=V_{n}^{*}\subset
z^{-n}\prod_{i=1}^{g}(z-\mu_{i})^{-1}\mathcal{P}\Leftrightarrow
V_{n}\supset z^{n}\prod_{i=1}^{g}(z-\mu_{i})\mathcal{P},
\end{equation*}
as desired. This finishes the proof of the lemma.
\end{proof}

\begin{lemma}\label{le8.5}                  
The tau function $\tau(n,t)$ associated with the rational flag
$\mathcal{V}$ defined by \eqref{8.5}, is given by
\begin{multline}                        \label{8.10}
\tau(n,t)=\prod_{j=1}^{g}\mu_{j}^{-n}\emph{exp}\Big(-\sum_{i=1}^{\infty}t_{i}
\sum_{j=1}^{g}\lambda_{j}^{i}\Big)\times\\
\emph{det}\Big(\lambda_{j}^{n+i-1}\emph{exp}(t,\lambda_{j})
-\delta_{j}\mu_{j}^{n+i-1}
\emph{exp}(t,\mu_{j})\Big)_{1\leq i,j\leq g}.
\end{multline}
\end{lemma}

\begin{proof} It follows easily from the definition of $V_{n}$ in \eqref{8.5}
that the wave function associated with this plane is
\begin{multline}                        \label{8.11}
w(n,t,z)=\frac{\mbox{exp}(t,z)}{\prod_{j=1}^{g}(z-\lambda_{j})}\times\\
\frac{\mbox{det}\big(\lambda_{j}^{n+i-1}\mbox{exp}(t,\lambda_{j})-
\delta_{j}\mu_{j}^{n+i-1}\mbox{exp}(t,\mu_{j});z^{n+i-1}\big)}
{\mbox{det}\big(\lambda_{j}^{n+i-1}
\mbox{exp}(t,\lambda_{j})-\delta_{j}\mu_{j}^{n+i-1}
\mbox{exp}(t,\mu_{j})\big)_{1\leq i,j\leq g}}.
\end{multline}
The notation in the numerator means the determinant of the
$(g+1)\times (g+1)$ matrix with $(i,j)$ entries as indicated there
for $1\leq i\leq g+1, 1\leq j\leq g$, and entries $z^{n+i-1},
1\leq i\leq g+1$, in the last column. Elementary row manipulations
with this determinant lead to
\begin{equation*}
w(n,t,z)=\frac{\mbox{exp}(t,z)}{\prod_{j=1}^{g}(z-\lambda_{j})}
z^{n+g}\frac{\prod_{j=1}^{g}(z-\lambda_{j})\tau(n,t-[z^{-1}])}{z^{g}\tau(n,t)},
\end{equation*}
in agreement with \eqref{8.3}, which defines the tau function. This
finishes the proof.
\end{proof}
\begin{remark}\label{re8.6}                     
The motivation for the (irrelevant) factor
$\prod_{j=1}^{g}\mu_{j}^{-n}$ in the definition \eqref{8.10} of
$\tau(n,t)$ will appear in the next subsection, see \eqref{8.36}.
\end{remark}

Having in view \thref{th8.1}, we are particularly interested in
rational curves with double points, for which the associated
commutative rings of difference operators contain a tridiagonal
operator.
\begin{proposition}\label{pr8.7}            
Let us assume that $\mu_{i}=\lambda_{i}^{-1}$,
and that the 2g points $\lambda_{i}, \mu_{i}, 1\leq i\leq g$, are
still distinct of each other, i.e. $\lambda_{i}\neq \pm 1$ and
$\lambda_{i}\neq \lambda_{j}^{\pm1},\forall\;i\neq j$. Then, the rank 1
commutative ring of difference operators $A_{\mathcal{V}}$, with
$V_{n}$ as in \eqref{8.5},
contains a tridiagonal operator $\hat{L}$ with one diagonal above and one
diagonal below the main diagonal. The spectral curve of the
ring has for equation
\begin{equation}                        \label{8.12}
\emph{Spec}(A_{\mathcal{V}}):\quad v^{2}=(u^{2}-4)\prod_{i=1}^{g}
\big(u-(\lambda_{i}+\lambda_{i}^{-1})\big)^{2},
\end{equation}
and the operator $\hat{L}$ can be obtained by iterating the
Darboux process starting with $L=E+E^{-1}$ as in \eqref{8.1}, and
adding bound states at the $g$ points $\lambda_{i}+\lambda_{i}^{-1}$.
\end{proposition}

\begin{proof} Since $\mu_{i}=\lambda_{i}^{-1}$, it is clear that
the function
\begin{equation}                        \label{8.13}
u=z+z^{-1},
\end{equation}
belongs to the ring $A_{\mathcal{V}}$ as defined in \eqref{8.2}, with
$uV_{n}\subset V_{n-1}$. This function has a simple pole both at
$z=0$ and at $z=\infty$ hence, by \thref{th8.3}, there exists a tridiagonal
operator $L_{u}$ with one diagonal above and one diagonal below
the main diagonal, satisfying $L_{u}w(n,t,z)=uw(n,t,z)$.

Since $u$ has a simple pole at $z=0$, all other generators of
$A_{\mathcal{V}}$ can be taken to be polynomials. The functions
\begin{equation}                        \label{8.14}
f_{k}=z^{-k}\prod_{i=1}^{g}(z-\lambda_{i})(z-\lambda_{i}^{-1}),\quad
k\in \mathbb{Z},
\end{equation}
which vanish at $\lambda_{i}$ and $\mu_{i}=\lambda_{i}^{-1}$,
obviously belong to $A_{\mathcal{V}}$. For $1\leq k\leq
g-1$, by subtracting an appropriate polynomial of degree $k$ in
$u$, one obtains polynomials (in $z$)
$q_{2g-1}(z),q_{2g-2}(z),\dots,q_{g+1}(z)$, of degrees $2g-1,
2g-2,\dots,g+1$, that belong to $A_{\mathcal{V}}$. Combining them
with the functions $f_{k},k\leq 0$, defined in \eqref{8.14}, we
conclude that $A_{\mathcal{V}}$ contains polynomials $q_{k}(z)$ of
degree $k$, $\forall k\geq g+1$. On the other hand, any
polynomial $q(z)\in A_{\mathcal{V}}$ of degree $k\leq g$, must be
identically constant. Indeed, such a polynomial must satisfy
$q(\lambda_{i})=q(\lambda_{i}^{-1}),\forall\;1\leq i\leq g$, which
amounts to a linear homogeneous system $Ac=0$, for the unknown
coefficients $c_{1},\dots,c_{g}$ of $q(z)=c_{0}+c_{1}z+\cdots+c_{g}z^{g}$.
One shows easily that
\begin{equation*}
\mbox{det}\;A=\prod_{i=1}^{g}\lambda_{i}^{-g}\prod_{i=1}^{g}(\lambda_{i}^{2}-1)
\prod_{1\leq i<j\leq g}(\lambda_{i}-\lambda_{j})\prod_{1\leq
i<j\leq g}(1-\lambda_{i}\lambda_{j}).
\end{equation*}
Since the $2g$ points $\lambda_{i},\mu_{i}=\lambda_{i}^{-1},1\leq
i\leq g$, are assumed to be distinct, this determinant is nonzero,
implying that $c_{1}=c_{2}=\cdots=c_{g}=0$.

The upshot of the discussion above is that the algebra
$A_{\mathcal{V}}$ is generated by $u$ in \eqref{8.13} and polynomials
of degree $k, k\geq g+1$, or equivalently by $u$ and the functions
$f_{k},k\in\mathbb{Z}$, introduced in \eqref{8.14}. We now show that
the functions $u$ and
\begin{equation*}
v=\frac{z-z^{-1}}{z^{g}}\prod_{i=1}^{g}(z-\lambda_{i})(z-\lambda_{i}^{-1}),
\end{equation*}
are enough to generate $A_{\mathcal{V}}$, which will establish
\eqref{8.12}. Clearly
\begin{equation*}
f_{g}=\prod_{i=1}^{g}\big(u-(\lambda_{i}+\lambda_{i}^{-1})\big)\quad
\mbox{and}\quad v=f_{g-1}-f_{g+1}.
\end{equation*}
Since $uf_{k}=f_{k-1}+f_{k+1}$, we deduce inductively that
$f_{g-1}, f_{g+1},f_{g-2},f_{g+2},\dots$, belong to the algebra
generated by $u$ and $v$.

It is well known that the spectral meaning of the curve \eqref{8.12} is
that the operator $L_{u}\equiv \hat{L}$, $u$ as in \eqref{8.13}, can be
obtained by iteration of the Darboux process, starting from $L=E+E^{-1}$ and
adding bound states at the points $\lambda_{i}+\lambda_{i}^{-1}$.
This can be easily deduced from the explicit form for the wave function
given in \eqref{8.11}, with $\mu_{i}=\lambda_{i}^{-1}$. The proof is complete.
\end{proof}

\subsection{The proof of \thref{th8.1}}
The proof of \thref{th8.1} is a combination of \prref{pr8.7}
above and \prref{pr8.8} below. It is convenient to introduce the
so-called reduced wave and adjoint wave functions, which are obtained
by omitting the factors $\mbox{exp}(t,z)$ and $\mbox{exp}^{-1}(t,z)$ in
\eqref{8.3} and \eqref{8.4} respectively:
\begin{equation}                        \label{8.15}
\Psi(n,z)=\mbox{exp}^{-1}(t,z)w(n,t,z),\quad
\Psi^{*}(n,z)=\mbox{exp}(t,z)w^{*}(n,t,z).
\end{equation}
We also omit to write the explicit dependence on $t$, which is
irrelevant for what follows; we just think of $t_{1},t_{2},t_{3},\dots$
as some free parameters.
\begin{proposition}\label{pr8.8}            
Consider a rational curve with double points,
obtained by identifying $2g$ distinct points of
$\mathbb{P}^{1}(\mathbb{C})$ in pairs
$\{\lambda_{i},\mu_{i}\},i=1,\dots,g$. Assume that
\begin{equation}                        \label{8.16}
\lambda_{i}=q^{k_{i}}\mu_{i},\;\emph{with}\;k_{i}\in\mathbb{N}\setminus\{0\}.
\end{equation}
Then the reduced wave function $\Psi(n,z)$, besides being a
common eigenfunction of a rank 1 commutative ring
$A_{\mathcal{V}}$ of difference operators (in $n$)
\begin{equation*}
L_{f}\Psi(n,z)=f(z)\Psi(n,z),\quad \forall f\in
A_{\mathcal{V}},
\end{equation*}
is also a common eigenfunction of a rank 1 commutative ring of
$q$-difference operators $A_{\tilde{\mathcal{V}}}$, in the variable z
\begin{equation}                        \label{8.17}
B_{g}\Psi(n,z)=g(n)\Psi(n,z),\quad \forall g\in
A_{\tilde{\mathcal{V}}}.
\end{equation}
\end{proposition}

Assuming for a moment this result we can prove \thref{th8.1}:
\begin{proof}[Proof of \thref{th8.1}]
From \eqref{8.16} we immediately obtain that the only intersection between
\prref{pr8.7} and \prref{pr8.8} is when
\begin{equation*}
\lambda_{i}=q^{k_{i}}\mu_{i}\;\mbox{and}\;\mu_{i}=\lambda_{i}^{-1}
\Leftrightarrow \lambda_{i}=\pm
q^{k_{i}/2},k_{i}\in\mathbb{N}\setminus\{0\},
\end{equation*}
in which case \eqref{8.12} reduces to \eqref{1.6}.
\end{proof}

As already announced, the proof of \prref{pr8.8} exploits an
interesting duality of the (reduced) wave function.
\begin{proof}[Proof of \prref{pr8.8}]
Using \eqref{8.16}, we have that
\begin{multline}                        \label{8.18}
\lambda_{j}^{n+i-1}\mbox{exp}(t,\lambda_{j})-\delta_{j}\mu_{j}^{n+i-1}
\mbox{exp}(t,\mu_{j})=\\
\mu_{j}^{n}\big\{x^{k_{j}}\lambda_{j}^{i-1}\mbox{exp}(t,\lambda_{j})-\delta_{j}
\mu_{j}^{i-1}\mbox{exp}(t,\mu_{j})\big\},\;\mbox{with}\;x=q^{n}.
\end{multline}
From \eqref{8.11} it follows that the (reduced) wave function is a
function of $x=q^{n}$ and $z$, i.e.
\begin{multline}                        \label{8.19}
\Psi(n,x)\equiv\Psi(x,z)=\\
\frac{e^{\frac{\log\;x\;\log\;z}{\log\;q}}}{\prod_{j=1}^{g}(z-\lambda_{j})}
\frac{\mbox{det}\big(x^{k_{j}}\lambda_{j}^{i-1}
\mbox{exp}(t,\lambda_{j})-\delta_{j}\mu_{j}^{i-1}
\mbox{exp}(t,\mu_{j});z^{i-1}\big)}{\mbox{det}\big(x^{k_{j}}\lambda_{j}^{i-1}
\mbox{exp}(t,\lambda_{j})-\delta_{j}\mu_{j}^{i-1}
\mbox{exp}(t,\mu_{j})\big)_{1\leq i,j\leq g}}.
\end{multline}

Any (reduced) wave function admits an expansion
\begin{equation}                        \label{8.20}
\Psi(x,z)=e^{\frac{\log\;x\;\log\;z}{\log\;q}}\Big(1+\frac{w_{1}(x)}{z}+
\frac{w_{2}(x)}{z^{2}}+\cdots\Big)\quad \mbox{as}\;z\to\infty.
\end{equation}
Since the $k_{i}$ are positive integers, it follows immediately
from \eqref{8.19} that
\begin{equation}                        \label{8.21}
\lim_{x\to\infty}e^{-\frac{\log\;x\;\log\;z}{\log\;q}}\Psi(x,z)=
\frac{\mbox{det}(\lambda_{j}^{i-1};z^{i-1})}{\prod_{j=1}^{g}(z-\lambda_{j})
\mbox{det}(\lambda_{j}^{i-1})_{1\leq i,j\leq g}}=1,
\end{equation}
showing that, when the special conditions \eqref{8.16} are satisfied,
we also have
\begin{equation*}
\Psi(x,z)=e^{\frac{\log\;x\;\log\;z}{\log\;q}}
\Big(1+\frac{\tilde{w}_{1}(z)}{x}+
\frac{\tilde{w}_{2}(z)}{x^{2}}+\cdots\Big)\quad \mbox{as}\;x\to\infty.
\end{equation*}
This suggests that the function $\tilde{\Psi}(x,z)$ obtained by
exchanging the variables $x$ and $z$ in $\Psi(x,z)$
\begin{equation}                        \label{8.22}
\tilde{\Psi}(x,z)=\Psi(z,x),
\end{equation}
should be the (reduced) wave function associated with
\emph{another} rational flag
$\tilde{\mathcal{V}}:\cdots\subset
\tilde{V}_{n+1}\subset\tilde{V}_{n}\subset\tilde{V}_{n-1}\subset\cdots$,
in the sense of \deref{de8.2}.
The correctness of this assertion will be established in \leref{le8.9}
below. Assuming the result, we deduce from \thref{th8.3} that
\begin{equation*}
B_{g}\tilde{\Psi}(x,z)\equiv\sum_{\text{finitely many
$i$'s$\;\in\Z$}}D_{i}(x)\tilde{\Psi}(q^{i}x,z)=g(z)\tilde{\Psi}(x,z),\quad
\forall \;g(z)\in A_{\tilde{\mathcal{V}}},
\end{equation*}
which, remembering the definitions of $\Psi(x,z)$ and
$\tilde{\Psi}(x,z)$ in \eqref{8.19} and \eqref{8.22}, amounts to
\begin{equation*}
\sum_{\text{finitely many $i$'s\;$\in\Z$}}D_{i}(z)\Psi(n,q^{i}z)
=g(q^{n})\Psi(n,z),
\end{equation*}
which establishes \eqref{8.17}. This finishes the proof of \prref{pr8.8}.
\end{proof}

As mentioned in the proof of \prref{pr8.8} given above, to
complete the argument, we still need to establish that $\tilde{\Psi}(x,z)$
defined as in \eqref{8.22} is indeed the (reduced) wave function of some
rational flag. We need some preliminaries.

We introduce the following multiplicative groups of formal
pseudodifference operators
\begin{equation}                        \label{8.23}
\mathcal{W}=\Big\{1+\sum_{i,j=1}^{\infty}w_{ij}x^{-j}D_{x}^{-i}\Big\},
\quad \mathcal{W}^{*}=\Big\{1+\sum_{i,j=1}^{\infty}w_{ij}x^{-j}D_{x}^{i}\Big\},
\end{equation}
with $w_{ij}\in\mathbb{C}$, $D_{x}$ and $D_{x}^{-1}$ the forward
and backward $q$-shift operators, $D_{x}h(x)=h(qx)$ and
$D_{x}^{-1}h(x)=h(q^{-1}x)$.

We denote by $\fa$ the adjoint isomorphism
\begin{equation}                        \label{8.24}
\fa:\mathcal{W}\to\mathcal{W}^{*}:W\curvearrowright
\fa(W)=(W^{-1})^{*},
\end{equation}
where $^{*}$ denotes the adjoint operator,
$(D_{x}^{i})^{*}=D_{x}^{-i}$. We define an anti-isomorphism
$\fb:\mathcal{W}\to\mathcal{W}$ by
\begin{multline}                        \label{8.25}
\fb(x)=D_{x},\quad \fb(D_{x})=x,\quad\mbox{i.e.}\\
\fb\Big(1+\sum_{i,j=1}^{\infty}
w_{ij}x^{-j}D_{x}^{-i}\Big)=1+\sum_{i,j=1}^{\infty}w_{ij}x^{-i}D_{x}^{-j},
\end{multline}
as well as an anti-isomorphism
$\fb^{*}:\mathcal{W}^{*}\to\mathcal{W}^{*}$ by
\begin{multline}                        \label{8.26}
\fb^{*}(x)=D_{x}^{-1},\quad \fb^{*}(D_{x})=x^{-1},\quad\mbox{i.e.}\\
\fb^{*}\Big(1+\sum_{i,j=1}^{\infty}w_{ij}x^{-j}D_{x}^{i}\Big)=
1+\sum_{i,j=1}^{\infty}w_{ij}x^{-i}D_{x}^{j}.
\end{multline}
With these definitions, it is straightforward to check that
\begin{equation}                        \label{8.27}
\fb^{*}\fa=\fa\fb.
\end{equation}

We can write \eqref{8.20} as
\begin{multline}                        \label{8.28}
\Psi(x,z)=W(x)\;e^{\frac{\log\;x\;\log\;z}{\log\;q}},\\
W(x)=1+w_{1}(x)D_{x}^{-1}+w_{2}(x)D_{x}^{-2}+\cdots.
\end{multline}
The functions $w_{i}(x)$ are rational functions of $x$
which, remembering \eqref{8.21}, satisfy
$\lim_{x\to\infty}w_{i}(x)=0$. Thus
\begin{equation}                        \label{8.29}
W(x)=1+\sum_{i,j=1}^{\infty}w_{ij}x^{-j}D_{x}^{-i} \quad
\mbox{as}\;x\to\infty,
\end{equation}
which shows that $W(x)\in \mathcal{W}$ as defined in \eqref{8.23}.

It is a general fact from the theory of the discrete KP-hierarchy
(see \cite{HI2}, \cite{I1}) that the adjoint (reduced) wave function as
defined by \eqref{8.4} and \eqref{8.15} can be expressed in terms of $W(x)$
as follows
\begin{equation}                        \label{8.30}
\Psi^{*}(x,z)=\fa(W)(q^{-1}x)\;e^{-\frac{\log\;x\;\log\;z}{\log\;q}},
\end{equation}
with $\fa(W)$ defined as in \eqref{8.24}. In the same papers, one can
find a proof  that \eqref{8.28} and \eqref{8.30} imply the so-called
\emph{bilinear identities}
\begin{equation}                        \label{8.31}
B\big(\Psi(q^{i}x,z),\Psi^{*}(x,z)\big)=0,\quad\forall\;i\geq 0,
\end{equation}
where $B$ is the residue pairing introduced in \eqref{8.7}.

With these preliminaries, we can establish the last needed
lemma.
\begin{lemma}\label{le8.9}                  
With the assumptions of \prref{pr8.8}, the new functions
\begin{equation}                        \label{8.32}
\tilde{\Psi}(x,z)=\Psi(z,x)\quad \emph{and}\quad
\tilde{\Psi}^{*}(x,z)=\Psi^{*}(qz,q^{-1}x)q^{-1}xz^{-1},
\end{equation}
are respectively the (reduced) wave and adjoint wave functions
associated with a rational flag $\tilde{\mathcal{V}}$.
\end{lemma}

\begin{proof} In terms of the anti-isomorphism $\fb$ introduced in
\eqref{8.25}, $\tilde{\Psi}(x,z)$ as defined in \eqref{8.32} can be written
as follows
\begin{equation}                        \label{8.33}
\tilde{\Psi}(x,z)=\fb(W)(x)\;e^{\frac{\log\;x\;\log\;z}{\log\;q}},
\end{equation}
with $W(x)$ as in \eqref{8.29}. Using \eqref{8.30} and the definitions
\eqref{8.26} and \eqref{8.32} of $\fb^{*}$ and $\tilde{\Psi}^{*}(x,z)$
respectively, one checks easily that
\begin{equation*}
\tilde{\Psi}^{*}(x,z)=\fb^{*}\big(\fa(W)\big)(q^{-1}x)\;
e^{-\frac{\log\;x\;\log\;z}{\log\;q}},
\end{equation*}
from which, using \eqref{8.27}, we deduce that
\begin{equation}                        \label{8.34}
\tilde{\Psi}^{*}(x,z)=\fa\big(\fb(W)\big)(q^{-1}x)\;
e^{-\frac{\log\;x\;\log\;z}{\log\;q}}.
\end{equation}

Let us define
\begin{align*}
\tilde{V}_{n}&=\mbox{Span}\{\tilde{\Psi}(n,z),\tilde{\Psi}(n+1,z),
\tilde{\Psi}(n+2,z),\dots\}\;\mbox{and}\\
\tilde{V}_{n}^{*}&=\mbox{Span}\{\tilde{\Psi}^{*}(n,z),
\tilde{\Psi}^{*}(n-1,z),\tilde{\Psi}^{*}(n-2,z),\dots\},
\end{align*}
with $\tilde{\Psi}(n,z)$ and $\tilde{\Psi}^{*}(n,z)$ the functions
obtained by substituting $q^{n}$ for $x$ in $\tilde{\Psi}(x,z)$ and
$\tilde{\Psi}^{*}(x,z)$ respectively. As recalled above (see \eqref{8.28},
\eqref{8.30} and \eqref{8.31}, with $W$ replaced by $\fb(W)$), the equations
\eqref{8.33} and \eqref{8.34} imply that $\tilde{\Psi}(n,z)$ and
$\tilde{\Psi}^{*}(n,z)$ satisfy the bilinear identities
$B\big(\tilde{\Psi}(n,z),\tilde{\Psi}^{*}(m,z)\big)=0, \;\forall\;n\geq m$.
Thus
\begin{equation}                        \label{8.35}
\tilde{V}_{n}^{*}=\tilde{V}_{n}^{\perp},
\end{equation}
where $\tilde{V}_{n}^{\perp}$ denotes the orthogonal of
$\tilde{V}_{n}$ with respect to $B$.

Using \eqref{8.18}, we can write $\tau(n,t)$ in \eqref{8.10} as
\begin{multline}                        \label{8.36}
\tau(n,t)\equiv
\tau(x,t)=\mbox{exp}\Big(-\sum_{i=1}^{\infty}t_{i}
\sum_{j=1}^{g}\lambda_{j}^{i}\Big)\times\\
\mbox{det}\Big(x^{k_{j}}\lambda_{j}^{i-1}
\mbox{exp}(t,\lambda_{j})-\delta_{j}\mu_{j}^{i-1}
\mbox{exp}(t,\mu_{j})\Big)_{1\leq i,j\leq g},
\end{multline}
with $x=q^{n}$, and thus $\tau(x,t)$ is a polynomial in $x$. Equations
\eqref{8.3} and \eqref{8.4} combined with the definition \eqref{8.32} of
$\tilde{\Psi}(x,z)$ and $\tilde{\Psi}^{*}(x,z)$, show then that
\begin{equation*}
\tilde{V}_{n}\subset z^{n}\tau(z,t)^{-1}\mathcal{P}\quad
\mbox{and}\quad \tilde{V}_{n}^{*}\subset
z^{-n}\tau(qz,t)^{-1}\mathcal{P},
\end{equation*}
which, because of \eqref{8.35}, is equivalent to
\begin{equation*}
z^{n}\tau(qz,t)\mathcal{P}\subset \tilde{V}_{n}\subset
z^{n}\tau(z,t)^{-1}\mathcal{P}.
\end{equation*}
Clearly, $\tilde{V}_{n}$ contains elements of orders
$n,n+1,n+2,\dots$, and only those, meaning that the codimension
of $\tilde{V}_{n}$ in $z^{n}\tau(z,t)^{-1}\mathcal{P}$ is equal to
the degree (as a polynomial in $z$) of $\tau(z,t)$. Thus the flag
$\tilde{\mathcal{V}}:\cdots\subset \tilde{V}_{n+1}\subset\tilde{V}_{n}\subset
\tilde{V}_{n-1}\subset\cdots$ is a rational flag in the sense of
\deref{de8.2}, and $\tilde{\Psi}(n,z)$ and $\tilde{\Psi}^{*}(n,z)$
are the (reduced) wave and adjoint wave functions associated with this flag.
The lemma is established.
\end{proof}

\begin{remark}\label{re8.10}                    
\prref{pr8.8} has been obtained by F.W. Nijhoff and O.A. Chalykh \cite{NC}
in the special case when
$\lambda_{i}=q\mu_{i}$, by writing the (reduced) wave function as
\begin{equation*}
\Psi(x,z)=e^{\frac{\log\;x\;\log\;z}{\log\;q}}\;\mbox{det}\Big\{\Id-
(x\;\Id+X)^{-1}[X,Z]_{q^{-1}}(z\;\Id+Z)^{-1}\Big\},
\end{equation*}
with $X$ and $Z$ some $g\times g$ matrices, so that the $q$-commutator
$[X,Z]_{q^{-1}}$ is a rank 1 matrix. This argument is insufficient to deal
with the (non-generic) situation that concerns us in \thref{th8.1}.
They conjecture that all bispectral rank 1 commutative rings of
difference operators can be parametrized in this way. In the limit
$q=1$, corresponding to rank 1 bispectral rings of differential
operators, this was proved in a highly non-trivial paper by G. Wilson
\cite{Wi2}. P. Iliev \cite{I2} has obtained independently of \cite{NC}
similar formulas.
\end{remark}


\end{document}